# Mathematical modeling of lymphocytes selection in the germinal center

**Vuk Milisic · Gilles Wainrib**

**Abstract** Lymphocyte selection is a fundamental operation of adaptive immunity. In order to produce B-lymphocytes with a desired antigenic profile, a process of mutation-selection occurs in the germinal center, which is part of the lymph nodes. We introduce in this article a simplified mathematical model of this process, taking into account the main mechanisms of division, mutation and selection. This model is written as a non-linear, non-local, inhomogeneous second order partial differential equation, for which we develop a mathematical analysis in the case of piecewise-constant coefficients. We assess, mathematically and numerically, the performance of the biological function by evaluating the duration of this production process as a function of several parameters such as the mutation rate or the selection profile, in various asymptotic regimes.

## Contents

## 1 Introduction

Understanding the immune system is a key challenge in current biological, medical and pharmaceutical research, leading to revolutionary biomedical applications such as vaccination, immunotherapy or specific antibody production.

Adaptive immunity is responsible for the evolution of the antibody repertoire, a learning process necessary to fight new foreign pathogens. This adaptation relies on a Darwinian process of Division-Mutation-Selection (DMS) occurring in the germinal centers [1], where an explosion of the mutation rate associated with B-cells division, called somatic hypermutation [2], is observed, hence providing a unique example of an evolutionary process occurring within living organisms.

Although the general qualitative description of this process is well-established in the literature [3,4,5], the quantitative assessment of this DMS process has remained largely unexplored experimentally, in particular due to the difficulty to gather precise phylogenetic data of the B-cell repertoire during the various phases [6]. Recently, several biological studies [7,8,9] have provided new experimental insights about the microscopic features of B-cell dynamics in the germinal center. Moreover, several key questions remain highly debated, such as the recycling of selected B-cells [10,11] or the neutrality of the mutation process. Over the last few years, this system has been studied using relatively detailed computational

V. Milisic
Université Paris 13, Laboratoire Analyse, Géométrie et Applications CNRS UMR 7539 99 av. Jean-Baptiste Clément, 93430 Villetaneuse, FRANCE E-mail: milisic@math.unvi-paris13.fr

G. Wainrib
Ecole Normale Supérieure, Département d'Informatique (DATA), 45 rue d'Ulm, 75005 Paris FRANCE. E-mail: gilles.wainrib@ens.fr

models and numerical simulations [12,13], providing a method to investigate several hypothesis and phenomena such as B-cell migration in the lymph node or the impact of recycling.

The aim of this article is to introduce a simplified macroscopic mathematical model of this process, in order establish rigorous mathematical foundations and to investigate theoretically the impact of a few key parameters, such as the mutation rate or the selection profile, on the performance of the B-cell production, characterized by the duration of the process or the final quality of the repertoire.

Mathematical modeling of population dynamics can be approached either from a microscopic agent-based point of view, considering the behavior of many individuals, leading to a stochastic system in high dimension, or from a macroscopic point of view, where global quantities such as the number of individuals in a given state are considered, leading to partial differential equations (PDE) or integro-differential equations. In this article, we focus on the latter approach, and introduce a PDE model describing a population of B-cells subject to division, mutation and selection processes. In our framework, the division and mutation features give rise to classical linear diffusion terms, whereas selection introduces an inhomogeneous term. When a sufficient amount of B-cell with desired properties has been selected, the overall process shall terminate, which is modeled through a feedback term. As a result, the proposed model is a non-linear, non-local and inhomogeneous elliptic PDE.

Several mathematical models of evolutionary dynamics with mutation and selection have been previously studied, especially in population dynamics [14],[15],[16], ranging from the early works of [17] to advanced mathematical models of adaptive dynamics. The model we introduce in this article has the particularity to combine spatial inhomogeities (in the space of traits) with a non-linear global feedback, which give rises to a specific PDE, for which we establish general properties, as well as specific estimates describing the impact of relevant parameters. In particular, our aim is to understand how the interplay between the mutation rate and the selection function determine the characteristic time-scale of the B-cell production process.

We first start with existence and uniqueness results for solutions of our new model. This gives a precise meaning to the solution, its regularity under the more general hypotheses. From the mathematical point of view, the main difficulty is due to the stiff, non-linear and non-local source term. The non-local feature is expressed in time and in the trait space.

In a second step, we compute, for a birth rate that is piecewise constant, the production time, *i.e*, the time $t_{\varrho_0}$ for which $\varrho_0$, a threshold selected population, is reached. This is performed with respect to two parameters of the model : $\varepsilon$ the width of the selection window, and $\mu$ the mutation rate.

The results of this study lead to several conclusions :

- for a general traits-square-integrable initial population of $B$ cels,
  - if $\varepsilon$ is small enough, and $\mu$ is greater than a threshold value, then $t_{\varrho_0}$ behaves as $|\log \varepsilon|$, and is independent on $\mu$,
  - for a fixed $\varepsilon$, if $\mu$ is small we observe that $t_{\varrho_0}$ tends towards a constant that depends on the initial birth rate, $\varepsilon$ and initial conditions.
  - for a fixed $\varepsilon$, when the mutation rate becomes large, $t_{\varrho_0}$ tends to another finite value, still depending on the same parameters.
- for an initial datum which is a Dirac mass whose support is a specific trait $z$ outside of the selection window :
  - if the domain is unbounded, we show that when $\mu$ tends to either 0 or $\infty$, $t_{\varrho_0}$ blows up,
  - at the contrary, when the domain is finite, $t_{\varrho_0}$ stays bounded for $\mu$ growing large.

  Intuitively these results can be interpreted saying that the larger the size of the repertoire is the more likely there should be an optimal mutation rate in terms of efficiency.

This distinction between different initial conditions comes from a debate in the literature, opposing supporters of mono-clonal germinal centers at the beginning of the process [18] to authors trying to prove and measure oligo-clonal initial populations [19].

The article is organized as follows. In section 2, we define mathematically our model and discuss its motivations, assumptions and limitations. In section 3, we derive general results concerning existence and uniqueness of solutions, as well as quantitative properties of the solutions. To gain further understanding into the dynamical behavior of the system with respect to the data, we study in section 4 the asymptotic behavior of $t_{\varrho_0}$ when $\varepsilon$ becomes small, whereas in 5 we consider for a fixed $\varepsilon$, the asymptotic regimes when $\mu$ is either large or small. This section focusses as well on different types of initial $B$-cell population. Throughout the paper, theoretical results are also illustrated with numerical simulations.

## 2 Mathematical model

We consider the time evolution of a population of lymphocytes during the Division-Mutation-Selection process within the germinal center.

### 2.1 Biological background

In this first section, we provide an elementary summary of the relevant biological background to describe our modeling approach. For the interested reader, we refer to [20] for a more detailed account.

The immune response to an external pathogen (virus, bacteria, etc.) involves many different types of cells and employs various strategies to eliminate the pathogenic sources.

One of the most important way to fight pathogens relies on the bonding between antigens and antibodies, which triggers an efficient immune response, recruiting many other agents such as macrophages or T-lymphocytes. Antibodies are macromolecular compounds, made of peptidic chains, and whose purpose is to bond with antigens, which are complementary molecules presented at the surface of pathogens. The antibody-antigens (A-A) bonding can be thought as putting a key in a lock, and is characterized by the concept of affinity, which quantifies the likelihood that this bonding occurs for a specific A-A pair.

The production of antibodies in the immune response is ensured by the B-type lymphocyte, which is an immune cell able to produce a single specific antibody. Therefore, it is essential for B-cells to be able to learn how to produce antibodies with a high affinity with a new given antigen. Studying this learning phenomenon is precisely the purpose of the present article.

Once a new antigen, called the target, is detected by the immune system, it is captured by the follicular dendritic cells, and brought to the lymph nodes. A simplified vision of the process of affinity maturation in the germinal center, a specific part of the lymph node, can be summarized as follows:

1. An initial population of a few immature B-cells enters in the germinal center. The affinity between the initial antibodies carried by these B-cells is in general relatively low, although it is not known whether this initial choice of B-cells is generic, or is already somehow adapted to the target.
2. During the first three days, B-cells divide and the population increases, until it migrates to another part of the germinal center, called the light zone. Notice that this population of B-cells is only able to produce a few types of different antibodies, typically less than 10.
3. In the light zone, the B-cells are now subject to a full division-mutation-selection process, on which is the focus of our model.
   - Division - Mutation : At each division, a single B-cell produces two cells, one of which has undergone a significant amount of mutations in the part of its DNA responsible for the production of the antibody peptide sequence. This process is called somatic hypermutation, since the mutation rate is now increased to extremely high levels, several orders of magnitude higher than in normal cell divisions.
   - Selection : follicular dendritic cells are in the light zone and present at their surface the target antigen. B-cells are also presenting their antibody at their surface and wander in a seemingly random manner in the light zone. Then, the affinity between the presented antibody and the target antigen determines the probability that a bonding occurs and lasts long enough. In that case, the B-cell is selected and receives a signal which enables the cell to escape from the germinal center. Otherwise, B-cells with a low affinity, which were not able to receive such a signal, die.
4. After the selection, a B-cell can have several fates : either it transforms into a plasma cell, which is able to produce and release antibodies to fight the pathogen, or into a memory cell whose aim is to remember the antigen (therefore being able to produce quickly high affinity antibodies, in case it comes back later). A third possible fate which is also discussed in the literature [11] is the possibility for the selected B-cell to come back inside the germinal center and to be subject a second time to the DMS process.
5. At some point, this affinity maturation terminates, and the precise biological mechanisms responsible for the determination of the stopping time remain unclear. However, it seems reasonable to consider that the process would stop as soon as a sufficient quantity of selected B-cell has escaped the germinal center.

### 2.2 Assumptions of the model

Of course, the above description is only a simplified and partial overview of a process which is in reality more complex and involves many different cell types, in particular T-cells. However, we think that it is neither possible nor wishful to take into account all the details of real biological processes into a mathematical model. Therefore, after running numerical simulations of various models (agent-based, stochastic models, PDE models with many variables), and motivated by the idea of introducing a simple, yet non-trivial, macroscopic mathematical model of the evolution of a population of B-cells during DMS phase, we have identified and selected what we consider to be the key parts of this complex process.

**Space of traits :** We consider that each B-cell is characterized by a *trait* corresponding to a specific antibody sequence. Instead of a discrete space of traits, composed with strings of amino-acids, we view here the trait as an abstract property of the antibody and we therefore consider the space of traits to be made with real numbers, for instance the interval $[0, 1]$.

**Affinity and selection :** The target antigen is also characterized by a trait $x_0$ in the same space, and if $x$ denotes the trait associated with a B-cell, then we consider that the A-A affinity is given by a function $s(x) = F(x_0, x)$ : the higher is $s(x)$, the higher is the likelihood that the B-cell with trait $x$ binds to the target. Since the affinity with the target summarizes all the necessary information about a given B-cell to decide its fate, we think that it is reasonable to consider real trait $x$, directly translated into an affinity through the function $F$. However, this model does not address the difficult question of understanding how a small change in the DNA of the B-cell will result into a change in the peptide chain of the antibody and finally into a modification of the affinity.

**Mutations :** In our modeling approach, we consider that a mutation will change slightly the trait $x$ to $x + dx$ in a diffusion manner, and that the affinity will also change slightly through the function $F$. This model of mutation does not take into account the possibility for a small mutation of the DNA to produce a large change in affinity. One way to overcome this difficulty would be to consider a non-local mutation kernel instead of a diffusion, but it seems relatively uneasy to us to make precise and justified assumptions on such a kernel.

**Termination :** We assume that the termination of the affinity maturation process is regulated by a measure of the number of selected B-cells. More precisely, we assume that the division rate is a decreasing function of the number of selected B-cells. Therefore, if this number reaches a certain value, the birth rate becomes lower than the death rate and the population inside the germinal center must extinct.

### 2.3 Mathematical model

We are now able to define precisely the mathematical model we will consider in this article. First, we introduce the following notations:

1. Parameters:
   - cell division rate function $Q : \mathbb{R} \to \mathbb{R}_+$ decreasing from $Q_0$ to $Q_1$
   - cell death rate $d > 0$
   - affinity-dependent selection function $s(.)$, peaked around $x_0$, the target.
   - mutation rate $\mu > 0$, which may either be a constant or a function of the trait $\mu : \mathbb{R} \to \mathbb{R}_+^*$
2. Variables:
   - $n(t, x)$ is the quantity of lymphocytes with an trait $x$
   - $\varrho(t)$ is the quantity of selected lymphocytes at time $t$ and is given by:

$$\varrho(t) = \int_0^t \int_{\mathbb{R}} s(x) n(z, x) dx dz$$

The domains in the trait-space will be denoted $\Omega$ and could be practically seen as a distance to a specific target trait. We are now able to formulate the main evolution equation, for $x$ in the trait space $\Omega := (0, 1)$ and $t \geq 0$. it is a initial boundary value problem reading : find the function $n$ solving

$$\begin{cases} \partial_t n(t,x) = (Q(\varrho(t)) - d - s(x))n(t,x) + \operatorname{div}(\mu \nabla n(t,x)), & (t,x) \in \mathcal{O}_T := \mathbb{R}_+ \times \Omega \\ \mu \partial_{\mathfrak{n}} n(t, \cdot) = 0, & (t,x) \in \Sigma_T := \mathbb{R}_+ \times \partial\Omega \\ n(0,x) = n_I(x), & \{0\} \times \Omega \end{cases} \quad (1)$$

where the second line is the homogeneous Neumann boundary condition ($\partial_{\mathfrak{n}} n := \nabla n \cdot \mathfrak{n}$), and the third one is the setting of initial data at time $t = 0$.

## 3 General results

In this section, we establish general results of existence and uniqueness, as well as spectral decomposition results, concerning system 1 under various assumptions on the coefficients.

### 3.1 Existence and uniqueness

In this sub-section, we present general existence and uniqueness results concerning system (1) (Theorem 1 for a Lipschitz continuous $Q$ and 2 for a piecewise continuous one).

We consider the Banach space

$$X = C\big([0,T]; L^2(\Omega)\big), \quad \|m\|_X := \sup_{0\le t\le T} \|m(t)\|_{L^2(\Omega)},$$

for some $T$ chosen later. Moreover we set

$$V(\mathcal{O}_T) := L^\infty((0,T); L^2(\Omega)) \cap L^2((0,T) \times H^1(\Omega))$$

and we define the form

$$I(t_1, n, \eta) := \int_\Omega n(t_1, x)\ \eta(t_1, x)dx - \int_0^{t_1}\int_\Omega n(t,x)\partial_t\eta(t,x)dxdt \\ + \int_0^{t_1}\int_\Omega \mu\nabla n(t,x)\cdot\nabla\eta(t,x) + (s(x) + d - Q(\varrho))n(t,x)\ \eta(t,x)dx$$

**Definition 31** *We call a weak solution of problem (1) any solution $n \in V(\mathcal{O}_T)$ s.t.*

$$I(T, n, \eta) = 0,$$

*for every function $\eta \in H^1(\mathcal{O}_T)$ s.t. $\eta \equiv 0$ when $t = 0$. We say moreover that the solution is consistent with the initial condition if*

$$I(T, n, \eta) = \int_\Omega n_I(x)\eta(0,x)dx, \quad \forall\eta \in H^1(\mathcal{O}_T).$$

*the latter equation will be denoted as variational formulation associated to the problem (1).*

**Hypotheses 31** *On the data we make the following assumptions*

*(i) The initial data $n_I(x)$ belongs to $L^2(\Omega)$ and is non-negative.*
*(ii) The function $Q$ is a uniformly Lipschitz function of its argument $\varrho$, $Q \in W^{1,\infty}(\mathbb{R})$.*
*(iii) The selection function $s$ is a bounded non-negative function of $x$.*
*(iv) The mutation rate $\mu$ is a bounded positive function of $x$.*
*(v) The death rate $d$ is a non-negative constant.*

**Theorem 1** *Under hypotheses 31, there exists a unique positive weak solution $n \in V(\mathcal{O}_T)$ for any positive time $T$.*

*Proof* We suppose in a first step that the constant $d$ is strictly positive. Then we prove the existence using the Banach fixed point theorem. We define $A$ a closed subset of $X$

$$A = \{m \in X, m > 0, \|m\|_X \le C_\Phi\},$$

where $C_\Phi$ is defined so that $\|n_I\|^2_{L^2(\Omega)} + \|Q\|^2_{L^\infty} T C^2_\Phi/(4d) < C^2_\Phi$. We denote by $c_1 := \|Q\|^2_{L^\infty}/(4d)$ For each $m \in A$, let $n$ be the weak solution associated to the problem :

$$\begin{cases} \partial_t n(t,x) - \mathrm{div}(\mu\nabla n(t,x)) + (d + s(x))n(t,x) = Q(\varrho(t))m(t,x), & (t,x) \in \mathbb{R}_+ \times \Omega \\ \varrho(t) = \int_0^t\int_\Omega s(x)m(l,x)dxdlx, & \\ \partial_{\mathfrak{n}} n(t,x) = 0 & t > 0, x \in \partial\Omega, \\ n(t=0,x) = n_I(x) > 0, & t = 0, x \in \Omega \end{cases} \tag{2}$$

This system defines the operator $\Phi : m \longmapsto n$. We prove that it admits a unique point in $A$.

- Existence :
  The solution $n \in V(\mathcal{O}_T)$ exists uniquely by standard parabolic theory (Theorem 5.1 p. 170 chap. III [21]) for any given $m \in L^\infty(\mathcal{O}_T)$. The $L^\infty$ bound follows the same way by Theorem 7.1 p. 181 [21]. Moreover by Theorem 4.2 p. 160 we know that
$$\int_0^{T-h} \|n(t+h,\cdot)-n(t,\cdot)\|^2_{L^2(\Omega)} dt = o(h).$$
- Non-negativeness :
  We suppose that $m$ is a positive function in $V(\mathcal{O}_T) \cap L^\infty(\mathcal{O}_T)$. We follow results from p. 183 [21] and we test in the weak formulation by $n_h^- := \min(n_h, 0)$ where $n_h$ is the Steklov approximation of $n$ i.e.
$$n_h(t,x) = \frac{1}{h}\int_t^{t+h} n(\tau,x)\,d\tau, \quad \forall (t,x) \in (0, T-h)\times \Omega.$$
  Applying such a test function is possible since $n_h^-$ actually does belong to $H^1(\mathcal{O}_T)$. Passing to the limit wrt the small parameter $h$ then gives the identity
$$\begin{aligned}\frac{1}{2}\left[\int_\Omega (n^-)^2(\tau,x)dx\right]_{\tau=0}^{\tau=t_1} - \int_0^{t_1}\int_\Omega \mu \nabla n \nabla n^- + (s(x)+d) n n^- dxdt \\ = \int_0^{t_1}\int_\Omega Q m n^- dxdt.\end{aligned}$$
  But the support of $n^-$ is the set where $n \le 0$ thus one has due to the positivity of $m$ and $Q$ that
$$\frac{1}{2}\left[\int_\Omega (n^-)^2(\tau,x)dx\right]_{\tau=0}^{\tau=t_1} \le 0.$$
  which gives that
$$\int_\Omega (n^-)^2(t_1,x)dx \le 0, \quad \forall t_1 \le T.$$
- Stability of the operator $\Phi$ in $X$ :
  By the same technique as above, we test by $n_h$ and pass to the limit wrt $h$ (for a more detailed explanation see p. 141-142 [21]) in the weak formulation, which writes:
$$\frac{1}{2}\left[\int_\Omega n^2(\tau,x)dx\right]_{\tau=0}^{\tau=t_1} + \int_{\mathcal{O}_{t_1}} \{\mu|\nabla n|^2 + (d+s(x))n^2\} dxdt = \int_{\mathcal{O}_{t_1}} Q(\varrho(t))\, m\, n\, dxdt,$$
  and applying Cauchy-Schwarz and the Young inequalities on the right hand side gives
$$\frac{1}{2}\left[\int_\Omega n^2(\tau,x)dx\right]_{\tau=0}^{\tau=t_1} \le \frac{\|Q\|_{L^\infty}T}{4d}\|m\|_X^2 \le \frac{\|Q\|_{L^\infty}T}{4d} C_\Phi^2, \quad \forall t_1 < T.$$
  Using the hypothesis on $C_\Phi$ we then deduce that in turn
$$\forall t \in [0,T], \|n(t,\cdot)\|_{L^2} \le C_\Phi,$$
  and thus $n \in A$.
- Second step: contraction. We denote $n_i = \Phi(m_i)$ for $i \in \{1,2\}$ where $m_i$ are two given functions in $A$. Then we denote $\tilde{n} := n_1 - n_2$ and $\tilde{m} := m_1 - m_2$, and we have
$$\begin{aligned}\frac{1}{2}\left[\int_\Omega \tilde{n}^2(\tau,x)dx\right]_{\tau=0}^{\tau=t_1} &+ \int_{\mathcal{O}_{t_1}} \{\mu|\nabla \tilde{n}|^2 + (d+s(x))\tilde{n}^2\} dxdt \\ &= \int_{\mathcal{O}_{t_1}} (Q(\varrho_1(t))m_1 - Q(\varrho_2)m_2)\tilde{n}\, dxdt, \\ &\le \frac{1}{4d}\int_{\mathcal{O}_{t_1}} \left[Q(\varrho_1)m_1 - Q(\varrho_2)m_2\right]^2 dxdt + d\int_{\mathcal{O}_{t_1}} \tilde{n}^2 dxdt\end{aligned}$$
  We estimate a bound for the right hand side
$$\frac{1}{4d}\int_{\mathcal{O}_{t_1}} \left[Q(\varrho_1)m_1 - Q(\varrho_2)m_2\right]^2 dxdt \le \frac{t_1\|Q\|^2_{L^\infty}}{2d}\|\tilde{m}\|_X^2 + \frac{\|m_2\|_X^2}{2d}\int_0^{t_1} |\tilde{Q}|^2 dt$$

The last term above can then be estimated using that

$$\int_0^{t_1} |\tilde{Q}|^2 dt = \int_0^{t_1} |Q(\varrho_1) - Q(\varrho_2)|^2 dt \le \|Q\|_{\text{Lip}}^2 \|s\|_{L^2(\Omega)}^2 \|\tilde{m}\|_X^2 \frac{t_1^3}{3}$$

which gives

$$\|\tilde{n}(t_1, \cdot)\|_{L^2(\Omega)}^2 \le t_1(c_1 + t_1^2 C_\Phi^2 c_2)\|\tilde{m}\|_X^2$$

where

$$c_2 = \frac{\|Q\|_{L^\infty}^2}{2d} \qquad c_3 = \frac{\|Q\|_{\text{Lip}}^2 \|s\|_{L^2(\Omega)}^2}{6d}$$

That finally provides

$$\|\tilde{n}\|_X \le \sqrt{T(c_2 + C_\Phi^2 T^2 c_3)}\|\tilde{m}\|_X$$

Choose $T$ small enough so that the contraction holds. The local existence and uniqueness follow from the Banach-Picard theorem.

– Third step: global existence
Denote $S_k = \sum_{i=1}^k T_i$. By induction, we assume that existence and uniqueness of (1) hold until the time $S_k$ with the corresponding bound $C_k$. Now consider the time interval $[S_k, S_{k+1}]$, the new problem is the variational formulation corresponding to the system written in a strong form :

$$\begin{cases} \partial_t \check{n}_{in}(t,x) = \big[Q(\varrho_{out}(t)) - d - s(x)\big]\check{n}_{in}(t,x) + \text{div}(\mu \nabla \check{n}_{in})(t,x), & (t,x) \in (S_k, S_{k+1}) \times \Omega, \\ \varrho_{out}(t) = \varrho_{out}(S_k) + \int_{[S_k,t]\times\Omega} s(x)\check{n}(l,x)dxdl, & \\ \check{n}_{in}(t = S_k, x) = n(S_k, x), & (t,x) \in (\{S_k\} \times \Omega, \\ \partial_{\mathrm{n}} \check{n}_{in}(t,0) = 0. & (t,x) \in ((S_k, S_{k+1}) \times \partial\Omega. \end{cases}$$

Iterating the same argument as in the first and second steps, existence and uniqueness hold on this new time interval if there exist $(T_k, C_k)_{k\in\mathbb{N}}$ s.t.

$$\begin{cases} \|n(S_k,\cdot)\|_{L^2}^2 + c_1 T_{k+1} C_{k+1}^2 < C_{k+1}^2, \\ \left(c_2 + C_{k+1}^2 T_{k+1}^2 c_3\right) T_{k+1} < 1, \end{cases} \qquad \forall k \in \mathbb{N}.$$

The first condition ensures that the operator $\Phi$ is an endomorphism while the second one insures its contractivity. Because $\|n(S_k,\cdot)\|_{L^2} \le C_k$, the first inequality holds if we suppose that

$$C_k^2 + c_1 T_{k+1} C_{k+1}^2 < C_{k+1}^2.$$

We choose $(T_k)$ and $(C_k)$ as

$$T_k := \frac{1}{2k\, c_1}, \quad C_k := \alpha k,$$

where $\alpha$ is s.t.

$$\alpha < \frac{c_1}{c_3}(4c_1 - c_2) = \frac{c_1}{2c_3 d}\|Q\|_{L^\infty}^2$$

the series $S_k = \sum_{i=1}^k T_i$ diverges. Thus global uniqueness and existence hold.

For a fixed positive constant $d$ we proved the theorem: there exists a unique weak solution $n_d \in V(\mathcal{O}_T)$ for any time $T$ solving

$$I_d(t_1, n_d, \eta) = 0, \quad \forall \eta \in H^1(\mathcal{O}_T)$$

vanishing on $t = 0$. Testing again with an appropriate averaged test function and passing to the limit gives :

$$\frac{1}{2}\left[\int_\Omega n_d^2(\tau,x)dx\right]_{\tau=0}^{\tau=t_1} + \int_{\mathcal{O}_{t_1}} \{\mu|\nabla n_d|^2 + s(x)n_d^2\}dxdt \le \int_{\mathcal{O}_{t_1}} Q(\varrho_d(t))\, n_d^2\, dxdt,$$

which provides by standard techniques

$$\|n_d\|_{V(\mathcal{O}_{t_1})} \le C e^{\|Q\|_{L^\infty} t_1} \|n_I\|_{L^2(\Omega)}$$

which is uniform wrt the size of $d$. By weak convergence arguments one passes to the limit when $d \to 0$. This proves the theorem in this specific case.

Hereafter we weaken the Lipschitz hypothesis made on $Q$ previously and we define

**Hypotheses 32** *We suppose that Hypotheses 31 (i),(iv) hold, moreover we suppose that*

*(ii)' $Q$ is a smooth function on $\mathbb{R}\setminus\{\varrho_0\}$, it admits two possibly different limits in the neighborhood of $\varrho_0$ $Q_\pm := \lim_{\varrho\to\varrho_0^\pm} Q(\varrho)$.*

*(iii)' $s$ is a positive definite function : $\inf_{x\in\Omega} s(x) > 0$.*

**Theorem 2** *Under Hypotheses 32, one has the same conclusions as in Theorem 1, except for uniqueness that only holds until a time $t_0^-$ defined below.*

*Proof* We define $Q_\delta$ a regularized non-linear source term

$$Q_\delta(\varrho) := \begin{cases} Q(\varrho) & \text{if } \varrho \in \mathbb{R}\backslash]\varrho_0, \varrho_0+\delta[, \\ Q_- + \frac{Q(\varrho_0+\delta)-Q_-}{\delta}(\varrho-\varrho_0) & \text{otherwise .}\end{cases}$$

As $Q_\delta$ is now Lipschitz continuous, one applies Theorem 1, denoting $n_\delta$ the corresponding unique solution, it belongs to $V(\mathcal{O}_T)$ for any time $T$, uniformly wrt $\delta$. Moreover, one has also that $n_\delta \in C(0,T;L^2(\Omega))$ uniformly wrt $\delta$. By uniformity of those bounds, one can extract a subsequence that we denote again $n_\delta$ that converges weakly in $L^2(0,T;H^1(\Omega))$ and weak star in $L^\infty(0,T;L^2(\Omega))$. For every fixed $t \geq 0$, $n_\delta(t,x)$ tends to $n(t,x)$ in $L^2(\Omega)$ weak. This shows that as $\delta \to 0$,

$$\begin{aligned}
&\int_\Omega n_\delta(t_1,x)\eta(t_1,x)dx \to \int_\Omega n(t_1,x)\eta(t_1,x)dx,\\
&\int_{\mathcal{O}_{t_1}} n_\delta(t,x)\partial_t\eta(t,x)dxdt \to \int_{\mathcal{O}_{t_1}} n(t,x)\partial_t\eta(t,x)dxdt,\\
&\int_{\mathcal{O}_{t_1}} \mu\nabla n_\delta(t,x)\nabla\eta(t,x)dxdt \to \int_{\mathcal{O}_{t_1}} \mu\nabla n(t,x)\nabla\eta(t,x)dxdt,\\
&\int_{\mathcal{O}_{t_1}} (s(x)+d)n_\delta(t,x)\eta(t,x)dxdt \to \int_{\mathcal{O}_{t_1}} (s(x)+d)n(t,x)\eta(t,x)dxdt.
\end{aligned}$$

It remains to show that

$$\int_{\mathcal{O}_{t_1}} Q_\delta(\varrho_\delta)n_\delta(t,x)\eta(t,x)dxdt \to \int_{\mathcal{O}_{t_1}} Q(\varrho)n(t,x)\eta(t,x)dxdt.$$

- $\varrho_\delta$ converges strongly to $\varrho$ in $C(0,T)$.
  Indeed, because $s\in L^\infty(\Omega)$,

  $$\rho_\delta(t) = \int_0^t\int_\Omega s(x)n_\delta(t,x)dxdt \to \int_0^t\int_\Omega s(x)n(t,x)dxdt =: \rho(t).$$

  and because $s$ and $n$ are non-negative functions $\varrho$ is a non-decreasing function.
- $\varrho_\delta$ is a strictly increasing function
  In the weak formulation $I(t,n_\delta,\eta)=0$, we set $\eta$ as a time dependent function independent on $x$ that solves

  $$\begin{cases} \partial_t\eta = (\|s\|_{L^\infty(\Omega)} + d - Q(\varrho_\delta(t)))\eta, & t\in(0,t_1)\\ \eta(0)=1, & t=0\end{cases}$$

  $\eta$ is then explicit and reads:

  $$\eta = \exp\left(\int_0^t (\|s\|_{L^\infty(\Omega)} + d - Q(\varrho_\delta(\tilde{t})))d\tilde{t}\right).$$

  We denote hereafter $\overline{n}$ the average of $n$ wrt to the trait variable :

  $$\overline{n} := \int_\Omega n(t,x)dx.$$

  Finally one has that:

  $$\overline{n}_\delta(t)\exp\left(\int_0^t (\|s\|_{L^\infty(\Omega)} + d - Q(\varrho_\delta(\tilde{t})))d\tilde{t}\right) = \overline{n}_\delta(0) + \int_{\mathcal{O}_t} (\|s\|_{L^\infty} - s(x))n_\delta(\tilde{t},x)\eta(\tilde{t})dx\,d\tilde{t}$$

the latter term being non-negative, and because the initial condition does not depend on $\delta$ one concludes that in fact:

$$\begin{aligned}\overline{n}_\delta(t) &\geq \overline{n}(0)\exp\left(\int_0^t (Q(\varrho_\delta(\tilde{t}) - (\|s\|_{L^\infty(\Omega)} + d))d\tilde{t}\right)\\ &\geq \exp\left(\left(\inf_{\varrho\in\mathbb{R}_+} Q(\varrho) - (\|s\|_{L^\infty} + d)\right)t\right) =: b(t) > 0\end{aligned}$$

which is positive definite for every finite time and the bound from below is uniform wrt $\delta$. Then one remarks that

$$\partial_t \varrho_\delta \geq \inf_{x\in\Omega} s(x)\overline{n}_\delta(t) \geq \inf_{x\in\Omega} s(x)b(t) > 0$$

which proves that $\varrho_\delta$ is strictly increasing for any fixed time uniformly wrt $\delta$ provided condition (iii)' of Hypotheses 32.

– Reaching $\varrho_0$ :
  As $\varrho_\delta$ is an increasing function whose initial datum is zero, there are two possibilities
  1. either $\varrho_\delta$ never reaches $\varrho_0$ *i.e.* $\varrho(t) < \varrho_0$ for every non negative $t \in \mathbb{R}$. Then for all times $Q_\delta(\varrho_\delta) = Q(\varrho_\delta)$ and $Q$ is always regular thus uniqueness results from Theorem 1 imply that $n_\delta = n$ a.e. in $(0,\infty)\times\Omega$. There is nothing to prove
  2. or there exists a time $t_0$ s.t. $\varrho_\delta(t_0) = \varrho_0$. Again by uniqueness, one has that

$$n_\delta(t,x) = n(t,x) \qquad \text{a.e } (t,x) \in (0,t_0)\times\Omega.$$

and thus this time $t_0$ is equal for every $\delta$. We fix a time $t_1 = 2t_0$ and define $\underline{c}$ as

$$\underline{c} := \inf_{t\in[0,t_1]} \overline{n}(0)\exp\left(\left(\inf_{\varrho\in\mathbb{R}_+} Q(\varrho) - (\|s\|_{L^\infty} + d)\right)t\right) > 0$$

We deduce then that in the neighborhood of $t_0$ one can write

$$\varrho_\delta(t_0+\omega) - \varrho_0 > \underline{c}\omega, \quad \forall\omega > 0.$$

Choosing then $\omega = \delta/\underline{c}$ one writes

$$\begin{aligned}\int_{\mathcal{O}_{t_1}} (Q_\delta(\varrho_\delta)n_\delta - Q(\varrho)n)\eta dxdt &= \int_{((0,t_0)\cup(t_0+\omega,t_1))\times\Omega} (Q(\varrho_\delta)n_\delta - Q(\varrho)n)\eta dxdt\\ &+ \int_{((t_0,t_0+\omega))\times\Omega} (Q_\delta(\varrho_\delta)n_\delta - Q(\varrho)n)\eta dxdt =: R_1 + R_2\end{aligned}$$

As previously shown $n_\delta(t,x) \equiv n(t,x)$ everywhere on $(0,t_0)$ a.e. $x\in\Omega$, thus $R_1$ reduces to

$$\begin{aligned}R_1 &= \int_{((t_0+\omega,t_1))\times\Omega} (Q(\varrho_\delta)n_\delta - Q(\varrho)n)\eta dxdt\\ &= \int_{((t_0+\omega,t_1))\times\Omega} (Q(\varrho_\delta)n_\delta - Q(\varrho)n_\delta + Q(\varrho)n_\delta - Q(\varrho)n)\eta dxdt\\ &\leq c_1\|\varrho_\delta - \varrho\|_{L^\infty(0,t_1)}\|n_\delta\|_{L^\infty((0,t_1);L^2(\Omega))} + \left|\int_{((t_0+\omega,t_1))\times\Omega}(n_\delta - n)Q(\varrho)\eta dxdt\right|\end{aligned}$$

where $c_1 = \|Q\|_{W^{1,\infty}(\mathbb{R}\setminus B(\varrho_0,\delta))}$, the latter term tends to zero when $\delta\to 0$ thanks to weak convergence arguments on $n_\delta$. On the other hand,

$$R_2 \leq \sqrt{\omega}\|Q\|_{L^\infty}\left(\|n\|_{L^\infty(0,t_1;L^2(\Omega))} + \|n_\delta\|_{L^\infty(0,t_1;L^2(\Omega))}\right)\|\eta\|_{H^1(\mathcal{O}_{t_1})}$$

This proves existence of a weak solution $n \in V(\mathcal{O}_T)$ of (1). Uniqueness follows up to the time $t_0^-$.

We end up with a time asymptotic result.

**Proposition 1** *Under hypotheses 31 and supposing that*

$$d > \lim_{\varrho\to\infty} Q(\varrho),$$

*the solution of (1) provides a monotone increasing function* $\varrho(t) = \int_0^t\int_\Omega s(x)n(\tilde{t},x)dxd\tilde{t}$ *that satisfies*

$$\varrho^\infty := \lim_{t\to\infty}\varrho(t) < \infty$$

*Proof* One has
$$\frac{d}{dt}\varrho = \int_{\mathbb{R}} s(x)n(t,x)dx \geq 0,$$
and thus $\varrho(t)$ is monotone increasing. By contradiction, assume that $\varrho(t) \to \infty$. Setting $\overline{n}(t) := \int_\Omega n(t,x)dt$ and testing the weak formulation with 1 gives

$$[\overline{n}(t) + \varrho(t)]_{t=t_1}^{t=t_2} = \int_{t_1}^{t_2} \left[Q(\varrho(t)) - d\right]\overline{n}(t)dt$$

Since $d > \lim_{t\to\infty} Q(\varrho(t))$, for $t_1$ large enough, the right hand side becomes negative implying $\overline{n}(t)+\varrho(t)$ shall not increase any more. This contradicts the assumption $\varrho(t) \to \infty$. We conclude that $\varrho^\infty < \infty$.

## 3.2 Spectral analysis

In this sub-section, we investigate the spectral decomposition of the solution, in particular under the assumption of piecewise-linear coefficient $Q(.)$.

### *3.2.1 Spectral decomposition*

As it is usual in the field of parabolic equations [22] one shall try the variable separation. To this aim we consider the spectral problem in the trait-space $\Omega$.

We denote $\mathcal{A}$ the operator defined by:

$$\mathcal{A}n(x) = -\mu\partial_{xx}^2 n(x) + s(x)n(x).$$

**Lemma 1** *We suppose that $s$ is a bounded function for almost every $x \in \Omega$. We denote the eigenvalues $\Lambda_k$ (resp. the eigenvectors $V_k$) of $\mathcal{A}$ the solution of the following equation:*

$$\begin{cases} \mathcal{A}V_k = \Lambda_k V_k, & a.e.\ x \in \Omega, \quad \forall k \in \mathbb{N} \\ \partial_x V_k(-1) = \partial_x V_k(1) = 0, & x \in \partial\Omega \end{cases} \tag{3}$$

*All eigenvalues of the spectrum are simple. The sequence of eigenvalues $(\Lambda_k)_{k\in\mathbb{N}}$ is monotone increasing and positive. The limit of $\Lambda_k$ when $k$ goes to infinity is infinite. Moreover the sequence $(V_k)_{k\in\mathbb{N}}$, it is an orthonormal basis of $L^2(\Omega)$. One can bound the eigenvalues : setting $m := \inf_\Omega \min(\mu, s)$ and $M := \sup_{x\in\Omega} \max(\mu, s)$*

$$m\left(\left(\frac{k\pi}{2}\right)^2 + 1\right) \leq \Lambda_k \leq M\left(\left(\frac{k\pi}{2}\right)^2 + 1\right).$$

*The first eigenvector is positive. The eigenvector associated to the eigenvalue $\Lambda_k$ has precisely $k$ zeros on $\Omega$.*

*Proof* The proof is standard [23,22] and can be found for instance in [24] Theorem 4.6.2 p. 87.

Then we project the initial data on the eigen-basis. We seek then the solution $n$ of system (1) as a time dependent superposition of modes :

$$n = \sum_{k\in\mathbb{N}} \alpha_k(t)V_k(x), \quad \forall(t,x) \in \mathbb{R}_+ \times \Omega$$

where the $\alpha_k$ should satisfy

$$\begin{cases} \partial_t\alpha_k + (\Lambda_k + d)\alpha_k = Q\left(\sum_k s_k \int_0^t \alpha_k(\tilde{t})d\tilde{t}\right)\alpha_k, \quad k \geq 1 \\ s_k = \int_\Omega V_k(x)s(x)dx = \Lambda_k \int_\Omega V_k dx, \\ \alpha_k(0) = \int_\Omega V_k(x)n_I(x)dx \end{cases} \tag{4}$$

We now make the assumption that the non-linearity $Q$ is piecewise-constant:

**Assumption 31** *The function* $Q(.)$ *is piecewise-constant :* $Q(\varrho) = Q_0$ *for* $\varrho \leq \varrho_0$ *and* $Q(\varrho) = Q_1$ *for* $\varrho > \varrho_0$.

Under this assumption, and for a time small enough $t < t_0$, *i.e* before $\rho(t)$ reaches $\varrho_0$, one solves explicitly (4) which gives:

$$\begin{cases} \alpha_k = \alpha_k(0)\exp\left((Q_0 - d - \Lambda_k)t\right), & t \in [0, t_0[ \\ \alpha_k(0) =< V_k, n_I > \end{cases}$$

where $< ., . >$ denotes the usual scalar product on $L^2(\Omega)$. We denote by $\Lambda_k^0 := b - \Lambda_k$ where $b := Q_0 - d$, so that we can write the solution as:

$$n(t,x) = \sum_{k=0}^{\infty} < n_I, V_k > e^{\Lambda_k^0 t} V_k(x), \quad t \leq t_0 \tag{5}$$

Let us define:

$$t_0 = \inf\{t \geq 0;\ \varrho(t) = \varrho_0\} \tag{6}$$

Then, one needs to estimate the value of $t_0$, from the above expression, and solve the system for time greater than $t_0$. To find the value of $t_0$, one needs to solve the non-linear equation $\varrho_0 = \varrho(t_0)$, more precisely:

$$\varrho_0 = \sum_k \phi_k \left(e^{\Lambda_k^0 t_0} - 1\right) \tag{7}$$

with

$$\phi_k := \frac{< n_I, V_k >< V_k, s >}{\Lambda_k^0}.$$

It appears relatively difficult to solve analytically this non-linear equation and find a formula for $t_0$. However, we will study in Sections 4 and 5 several approaches to derive asymptotic estimates of this key quantity.

It remains to solve the system for time $t$ larger than $T$. The strategy is exactly the same, except that the new eigenvalues are now $\Lambda_k^1 = (Q_1 - d) - \Lambda_k$, and that the new initial condition is $n(T, x)$, denoted $n_T(x)$. Therefore, for $t > T$:

$$n(t,x) = \sum_k < n_T, V_k > e^{\Lambda_k^1 (t-T)} V_k(x) \tag{8}$$

Since $(\Lambda_k)_{k\in\mathbb{N}}$, the spectrum of $\mathcal{A}$ is located on the positive real axis. Therefore, we deduce the following time-asymptotic result: if $d > Q_1$ then

$$\lim_{t\to\infty} n(t,x) = 0, \quad a.e. x \in \Omega$$

*3.2.2 Spectral calculus*

We now consider the problem of finding explicit expressions for the eigenvalues and eigenvectors associated with problem (3).

To this end, we now make some more assumptions on $Q$, $s$ and $\mu$ :

**Assumption 32** *1. The function* $Q(.)$ *is piecewise-constant :* $Q(\varrho) = Q_0$ *for* $\varrho \leq \varrho_0$ *and* $Q(\varrho) = Q_1$ *for* $\varrho > \varrho_0$
*2. The function* $s(.)$ *is piecewise-constant :* $s(x) = s_\varepsilon(x) = 1$ *if* $x \in [-\varepsilon, \varepsilon]$ *and* $0$ *otherwise.*
*3. The function* $\mu(.)$ *is constant* $\mu(x) = \mu > 0$.

Notice that $\varepsilon$ is a parameter in $(0, 1]$ that we do not assume, here, to be small. We consider the limit of small $\varepsilon$ in Section 4.

We study the eigen-problem (3) under Assumptions 32. This means find the sequence $(V_k, \Lambda_k)$ solving

$$-\mu V_k'' + s_\varepsilon V_k = \Lambda_k V_k, \text{ in } ]-1,1[, \quad V_k'(-1) = V_k'(1) = 0.$$

Because of the symmetry of the problem we restrict the domain to $]0,1[$ and write

$$-\mu V_k'' + s_\varepsilon V_k = \Lambda_k V_k, \text{ in } ]0,1[, \quad V_k'(0) = V_k'(1) = 0.$$

On each part of the domain one has a constant coefficient problem that can be solved. Our goal is to construct by composition the complete eigenproblem in this particular case. Converting the second order problem into a first order system one has:

$$\begin{cases} \partial_x \mathbf{Y}(x,\Lambda) = M(x,\Lambda)\mathbf{Y}(x,\Lambda), & \text{in } ]0,1[ \\ Y_2(0,\Lambda) = Y_2(1,\Lambda) = 0, & x \in \{0,1\} \end{cases} \tag{9}$$

where $\mathbf{Y} = (V_k, \partial_x V_k)^T$. So the complete problem can be solved by piecewise exponentials:

$$\mathbf{Y}(x,\Lambda) = \Phi(x,\Lambda)\mathbf{Y}(0), \quad \Phi(x,\Lambda) := \begin{cases} e^{M(\Lambda-1)x} & \text{if } x < \varepsilon \\ e^{M(\Lambda)(x-\varepsilon)} e^{M(\Lambda-1)\varepsilon} & \text{if } x > \varepsilon \end{cases},$$

and

$$M(\Lambda) = \begin{pmatrix} 0 & 1 \\ -\omega_0^2(\Lambda) & 0 \end{pmatrix} \mathbb{1}_{(0,\varepsilon)}(x) + \begin{pmatrix} 0 & 1 \\ -\omega_1^2(\Lambda) & 0 \end{pmatrix} \mathbb{1}_{(\varepsilon,1)}(x)$$

where $\omega_0^2 := -\frac{s_0-\Lambda}{\mu}$ and $\omega_1^2 := \frac{\Lambda}{\mu}$. The exponential matrices are then explicit:

$$e^{M(\Lambda)x} = \begin{pmatrix} \cos(\omega_0 x) & \sin(\omega_0 x)/\omega_0 \\ -\omega_1 \sin(\omega_0 x) & \cos(\omega_0 x) \end{pmatrix} \mathbb{1}_{(0,\varepsilon)}(x) + \begin{pmatrix} \cos(\omega_1 x) & \sin(\omega_1 x)/\omega_1 \\ -\omega_1 \sin(\omega_1 x) & \cos(\omega_1 x) \end{pmatrix} \mathbb{1}_{(\varepsilon,1)}(x).$$

The boundary problem (9) can be rewritten in an algebraic form :

$$T\mathbf{Y}(0) = 0, \quad T := N_1 + N_2\Phi(1), \quad N_1 := \begin{pmatrix} 0 & 1 \\ 0 & 0 \end{pmatrix}, \quad N_2 := \begin{pmatrix} 0 & 0 \\ 0 & 1 \end{pmatrix}$$

that has a solution when $\det T = 0$ i.e.

$$\left| N_1 + N_2 e^{M(\Lambda)(1-\varepsilon)} e^{M(\Lambda-1)\varepsilon} \right| = 0$$

this leads to a simpler condition

$$\omega_0 \tan_0 + \omega_1 \tan_1 = 0$$

where $\tan_j := \tan(\omega_j x_j)$, $j \in \{0,1\}$, $x_0 = \varepsilon$ and $x_1 = (1-\varepsilon)$. As a function of $\Lambda$ the latter equation becomes:

$$\sqrt{\Lambda - 1} \tan\left(\delta\sqrt{\Lambda - 1}\mu\right) + \sqrt{\Lambda}\tan((1-\delta)\sqrt{\Lambda}) = 0. \tag{10}$$

The corresponding eigenvector is

$$V(x) = \mathbb{1}_{x\leq\varepsilon}\cos(x\omega_0) \quad + \mathbb{1}_{x\geq\varepsilon}\left\{ \cos(\varepsilon\omega_0)\cos(\omega_1(x-\varepsilon)) - \sin(\varepsilon\omega_0)\sin(\omega_1(x-\varepsilon))\frac{\omega_0}{\omega_1} \right\}, \tag{11}$$

and one should simply take into account whether $\omega_0$ is a pure imaginary or real number in order to pass from hyperbolic to standard trigonometric functions. So there exists a sequence $(\Lambda_k, V_k)_{k\in\mathbb{N}}$, s.t. for each $k \in \mathbb{N}$, $\Lambda_k$ solves (10) and $V_k$ writes as in (11). But, for a given set of data $s(x)$ and $\mu$, (10) the solution $\Lambda_k$ is not explicit. In what follows we approximate it by a spectral asymptotic expansion.

## 4 Asymptotic analysis for narrow selection profiles

In this section, we focus our attention on the case of a narrow selection profile, namely considering that the selection function $s(x)$ can be written as $s(x/\epsilon)$ where $\epsilon$ is a small parameter. In biological terms, this assumption means that selection is very specific, and that only B-cells with a trait very similar to the target are selected.

We first consider this asymptotic regime from the spectral point of view, and then construct an asymptotic expansion of the solution in $\epsilon$. These two approaches enable an asymptotic estimation of the time to threshold, characterizing the duration of the production process until a sufficient level of selected B-cell is reached.

## 4.1 Asymptotic spectral analysis

### *4.1.1 Spectrum of auto-adjoint operators with compact inverse*

Because the Neumann problem admits a zero eigenvalue we shift the spectrum by adding the identity because then $-\Delta + I$ has an auto-adjoint compact inverse and the spectral theory can be used. Since the shift is artificial and does not change any of the results presented below we omit it : in what follows we return to the original operator. Any time that we mention that the operator has a compact inverse, it is understand in the sense above (*i.e.* modulo a unit shift).

In a first step we study again problem (3) under assumptions 32. Since the operator $\mathcal{A}$ is auto-adjoint and has a compact inverse from $L^2(0,1)$ into itself $L^2(0,1)$ it admits (see Theorem 6. p. 38 [22]) a discrete spectrum that can be arranged into an increasing sequence of real eigenpairs denoted $(\Lambda_{\varepsilon,k}, V_{\varepsilon,k})$ for $k \in \mathbb{N}$.

### *4.1.2 An asymptotic Ansatz*

We develop $V_{\varepsilon,k}$ using an asymptotic expansion, to this purpose we approach $V_{\varepsilon,k}$ by a series reading

$$\mathcal{V}_{\varepsilon,k}(x,y) := v_{0,k}(x,y) + \varepsilon v_{1,k}(x,y) + \varepsilon^2 \pi_{\varepsilon,2,k}(x,y)$$

where $x$ represents the slow variable ($x \in (0,1)$) and $y$ the fast variable (typically $y = x/\varepsilon$). The asymptotic expansion provides a set of equations that finally provide that $v_{0,k}$ does not depend on $y$ and solves

$$\begin{cases} -\mu v''_{0,k} = \lambda_{0,k} v_{0,k}, & x \in (0,1) \\ v'_{0,k}(0) = v'_{0,k}(1) = 0 \end{cases} \tag{12}$$

the solution is explicit : $\lambda_{0,k} := \mu(k\pi)^2$ for $k \in \mathbb{Z}$, whereas normalizing the eigenvectors gives :

$$v_{0,k}(x) = \begin{cases} 1, & \text{if } k = 0 \\ \sqrt{2}\cos(k\pi x), & \text{otherwise} \end{cases}$$

Since $v_{0,k}$ does not solve the original equation, we correct it by adding a second order microscopic corrector $\pi_{\varepsilon,2,k}(x) := \pi_2(x/\varepsilon) v_{0,k}(x)$ where $\pi_2(y)$ solves

$$\begin{cases} -\mu \partial^2_{y^2} \pi_2 + s_0(y) = 0, & x \in (0,\infty) \\ \partial_y \pi_2(0) = 0. \end{cases} \tag{13}$$

The solution is explicit:

$$\pi_2(y) = \pi_2(0) + \frac{1}{\mu} \int_0^y \int_0^z s(\tilde{z}) d\tilde{z} dz.$$

Since all the functions are defined up to a constant we omit them in the rest of the section. If we take the first order derivative of $\pi_{2,k}$ one has:

$$\partial_y \pi_{2,k}(y) = \frac{1}{\mu} \begin{cases} \int_0^y s(z)dz & \text{if } y < 1 \\ \overline{s} := \int_0^1 s(z)dz & \text{otherwise} . \end{cases}$$

In order to reduce the contribution of the growth at infinity of the latter microscopic boundary layer function, we introduce the first order spectral problem: $v_{1,k}$ solves

$$\begin{cases} -\mu v''_{1,k} = \lambda_{0,k} v_{1,k} + \lambda_{1,k} v_{0,k} + 2\overline{s} v'_{0,k}, & x \in (0,1), \\ v'_{1,k}(0) = 0, & \text{if } x = 0, \\ v'_{1,k}(1) = -\dfrac{\overline{s}}{\mu} v_{0,k}(1), & \text{if } x = 1. \end{cases}$$

Because the operator admits a kernel of dimension one, the data of this problem should be polar to the kernel of the adjoint. Therefore, due to this latter condition:

$$\lambda_{1,k} := \frac{\overline{s}\, v^2_{0,k}(0)}{\|v_{0,k}\|^2_{L^2(0,1)}} = \begin{cases} 2\overline{s} & \text{if } k \neq 0 \\ \overline{s} & \text{otherwise} \end{cases}$$

Interestingly enough this first order eigen-contribution is independent on $\mu$. The solution of the latter problem then reads (modulo a multiple of $v_{0,k}$)

$$v_{1,k} = -\frac{\bar{s}}{\mu} \cdot \begin{cases} x^2/2 & \text{ăif ă} k = 0, \\ (1-x)\frac{v'_{0,k}(x)}{(k\pi)^2} + x v_{0,k}(x) & \text{otherwise} \end{cases}$$

We compute the problem solved by $\mathcal{V}_{\varepsilon,k}$ :

$$\begin{aligned} -\mu \mathcal{V}''_{\varepsilon,k} + s_\varepsilon \mathcal{V}_{\varepsilon,k} = & -\mu v''_{0,k} + s_\varepsilon v_{0,k} - \varepsilon \mu v''_{1,k} + \varepsilon s_\varepsilon v_{1,k} \\ & - \varepsilon^2 \mu \left( \frac{\partial^2_{y^2} \pi_{2,k}(x/\varepsilon)}{\varepsilon^2} v_{0,k} + 2 \frac{\partial_y \pi_{2,k}(x/\varepsilon)}{\varepsilon} v'_{0,k} + \pi_{2,k} v''_{0,k} \right) + \varepsilon^2 s_\varepsilon \pi_{\varepsilon,2,k} \\ = & (\lambda_{0,k} + \varepsilon \lambda_{1,k})(v_{0,k} + \varepsilon v_{1,k} + \varepsilon^2 \pi_{\varepsilon,2,k}) - \varepsilon^2 (\lambda_{1,k} v_{1,k} + \varepsilon \pi_{\varepsilon,2,k}) \\ & + \varepsilon s_\varepsilon (v_{1,k} + \varepsilon \pi_{\varepsilon,2,k}) + 2\varepsilon \left( \int_x^\varepsilon s(z/\varepsilon) dz \mathbb{1}_{(0,\varepsilon)}(x) \right) v'_{0,k} \\ = & \lambda_{\varepsilon,1,k} \mathcal{V}_{\varepsilon,k} + \varepsilon \left( s_\varepsilon v_{1,k} + 2 \left( \int_x^\varepsilon s(z/\varepsilon) dz \mathbb{1}_{(0,\varepsilon)}(x) \right) v'_{0,k} \right) + O(\varepsilon^2). \end{aligned}$$

where we defined $\lambda_{\varepsilon,1,k} := \lambda_{0,k} + \varepsilon \lambda_{1,k}$.

**Theorem 3** *Under the hypotheses on $s_\varepsilon$, the tuple $(\lambda_{\varepsilon,1,k}, \mathcal{V}_{\varepsilon,k})$ is a generalized eigen-pair,* i.e. *it verifies*

$$|a_\varepsilon(\mathcal{V}_{\varepsilon,k}, v) - \lambda_{\varepsilon,1,k}(\mathcal{V}_{\varepsilon,k}, v)| \leq C(k+1)\pi \varepsilon^{\frac{3}{2}} \|\mathcal{V}_{\varepsilon,k}\|_{L^2(0,1)} \|v\|_{L^2(0,1)}, \forall v \in H^1(0,1),$$

*which implies that*

$$|\Lambda_{\varepsilon,k} - \lambda_{\varepsilon,1,k}| \leq C(k+1)\pi \varepsilon^{\frac{3}{2}}.$$

*The constant $C$ depends on $s$, but not on $\mu$.*

*Proof* We use the problem solved by $\mathcal{V}_{\varepsilon,k}$ computed above and estimate the first order terms :

$$\begin{aligned} \varepsilon \left( \left( s_\varepsilon v_{1,k} + 2 \left( \int_x^\varepsilon s(z/\varepsilon) dz \mathbb{1}_{(0,\varepsilon)}(x) \right) v'_{0,k} \right), v \right)_{L^2(0,1)} \leq & \varepsilon \|s_\varepsilon\|_{L^2(0,1)} \|v_{1,k}\|_{L^\infty(0,1)} \|v\|_{L^2(0,1)} \\ & + 2\varepsilon \|s\|_{L^\infty} \|\mathbb{1}_{0,\varepsilon}\|_{L^2(0,1)} \left\| v'_{0,k} \right\|_{L^\infty(0,1)} \|v\|_{L^2(0,1)} \\ \leq & C \varepsilon^{\frac{3}{2}} (1+k)\pi \|\mathcal{V}_{\varepsilon,k}\|_{L^2(0,1)} \|v\|_{L^2(0,1)} \end{aligned}$$

the higher order terms are bounded in the $L^\infty(0,1)$ norm wrt $\varepsilon$ and independent on $k$, then applying Theorem 6 p. 38 in [22], whose proof by contradiction based on the spectral decomposition of $a_\varepsilon$ can be found in [25,26].

*4.1.3 Numerical comparison*

In fig. 1, we compare the first order approximation of the eigenvalues for the first ten values of $k$ and we display the numerical solution of the characteristic polynomial (10). We choose $\varepsilon = 0.1$.

*4.1.4 Application: estimation of the time to threshold*

We consider the non-linear equation:

$$\varrho_0 = \sum_k \Phi_{\varepsilon,k} \left( e^{(b-\Lambda_{\varepsilon,k})t_{\varrho_0}} - 1 \right) \tag{14}$$

with

$$\Phi_{\varepsilon,k} := \frac{< n_I, V_{\varepsilon,k} >< V_{\varepsilon,k}, s >}{b - \Lambda_{\varepsilon,k}}.$$

There exists no analytical solution for $t_{\varrho_0}$ for the above equation. However, it is possible to obtain an asymptotic description of $t_{\varrho_0}$. In Section 3.2, we derived an asymptotic expansion of the eigenvalues $\Lambda_{\varepsilon,k}$ and eigenvectors $V_{\varepsilon,k}$ when $\epsilon \to 0$ in the form:

$$\Lambda_{\varepsilon,k} = \sum_{l \geq 0} \lambda_{l,k} \epsilon^l \text{ and } V_{\varepsilon,k} = \sum_{l \geq 0} v_{l,k} \epsilon^l \tag{15}$$

We first show the following:

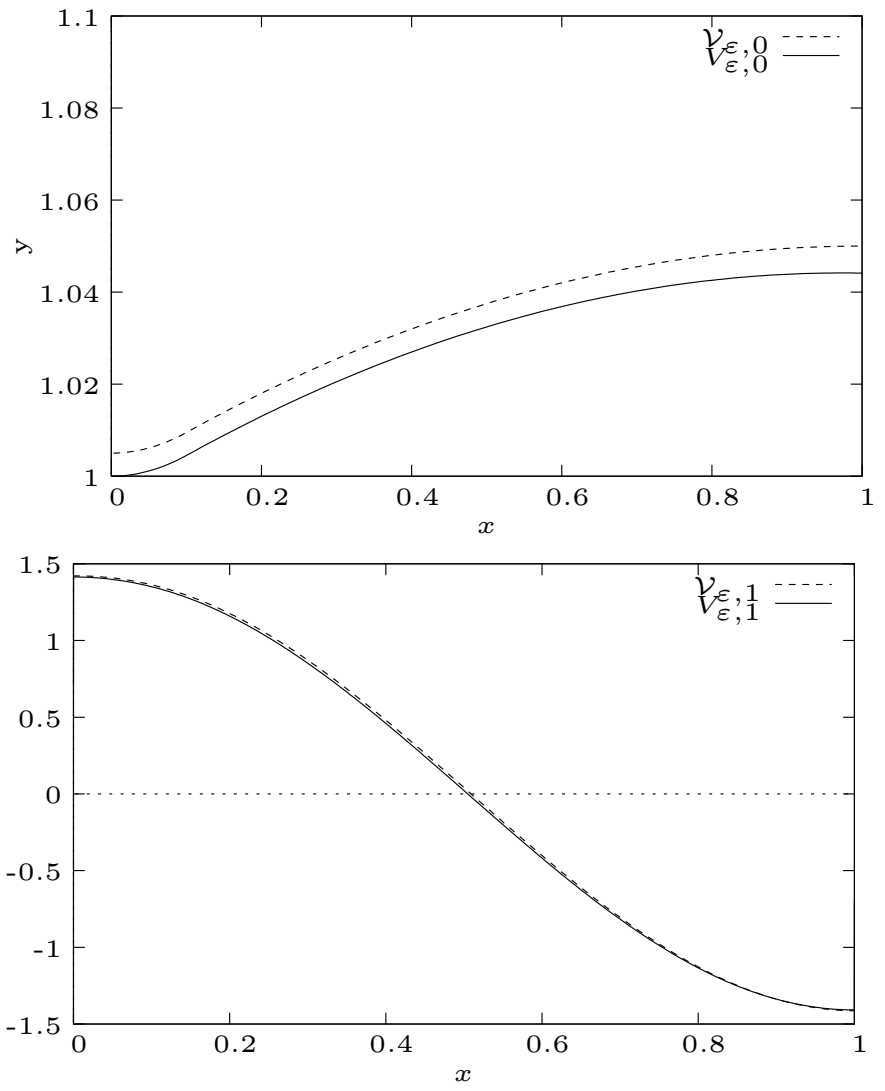

| $k$ | $\lambda_{\varepsilon,1,k}$ | $\Lambda_{\varepsilon,k}$ |
|---|---|---|
| 0 | 0.1 | 0.09735591831 |
| 1 | 10.0696044 | 10.06329604 |
| 2 | 39.6784176 | 39.6548511 |
| 3 | 89.02643961 | 88.97754925 |
| 4 | 158.1136704 | 158.0375159 |
| 5 | 246.94011 | 246.840386 |
| 6 | 355.5057584 | 355.3903082 |
| 7 | 483.8106157 | 483.6890491 |
| 8 | 631.8546817 | 631.7357679 |
| 9 | 799.6379565 | 799.5275421 |
| 10 | 987.1604401 | 987.0604181 |
| 11 | 1194.422133 | 1194.330628 |

**Fig. 1** Numerical comparison between asymptotic and numerical eigen-pairs. On the left the eigenvectors for $k \in \{0,1\}$ and $\varepsilon = 0.1$, on the right the eigenvalues for $k \in \{0,\dots,11\}$.

**Lemma 2** *Assuming that the following limits exists:*

$$\lim_{\varepsilon\to 0} \frac{1}{\varepsilon} \int s(x/\varepsilon) v_{i,k}(x)dx = \underline{v}_{i,k} \tag{16}$$

*and that* $\lambda_{0,k} \neq b$ *for all* $k \geq 0$*, then the coefficient* $\Phi_{\varepsilon,k}$ *can be approximated by:* $\Phi_{\varepsilon,k} = \varepsilon\phi_{0,k} + O(\varepsilon^2)$*, with* $\phi_{0,k} = \frac{<n_I, v_k^0> \underline{v}_k^0}{b-\lambda_{0,k}}$.

*Proof* By definition,

$$\Phi_{\varepsilon,k} := \frac{< n_I, v_k^\epsilon >< s^\epsilon, v_k^\epsilon >}{b - \Lambda_k^\epsilon} \tag{17}$$

From the asymptotic description of the spectrum, we know that:

$$\Lambda_{\varepsilon,k} = \lambda_{0,k} + \varepsilon\lambda_{1,k} + O(\varepsilon^2) \tag{18}$$

$$V_{\varepsilon,k} = v_{0,k} + \epsilon v_{1,k} + O(\varepsilon^2) \tag{19}$$

Therefore, assuming that the following limits exists:

$$\lim_{\epsilon\to 0} \frac{1}{\epsilon} \int s(x/\epsilon) v_{i,k}(x)dx = \underline{v}_{i,k} \tag{20}$$

we obtain:

$$< s_\epsilon, V_{\varepsilon,k} >= \varepsilon\underline{v}_{0,k} + \epsilon^2\underline{v}_{1,k} + O(\epsilon^3) \tag{21}$$

Gathering the above estimates, we obtain the desired result.

An asymptotic description of $t_{\varrho_0}$ is given by the following result, which essentially assumes that only the first mode is growing (assumption (iii) below).

**Lemma 3** *If* $t_{\varrho_0}$ *is the solution of equation* (14)*, with:*

*(i)* $\Lambda_{\varepsilon,k} = \lambda_{0,k} + O(\epsilon)$
*(ii)* $\Phi_{\epsilon,k} = \epsilon\phi_{0,k} + O(\epsilon^2)$
*(iii)* $b - \lambda_{0,0} > 0$ *and* $b - \lambda_{0,k} < 0$ *for all* $k > 0$

*Then* $t_{\varrho_0}$ *diverges to* $+\infty$ *when* $\epsilon \to 0$ *as:*

$$t_{\varrho_0} = \frac{1}{b - \lambda_{0,0}} \ln\left(\frac{\varrho_0}{\phi_{0,0}}\frac{1}{\epsilon} + o\left(\frac{1}{\epsilon}\right)\right) \tag{22}$$

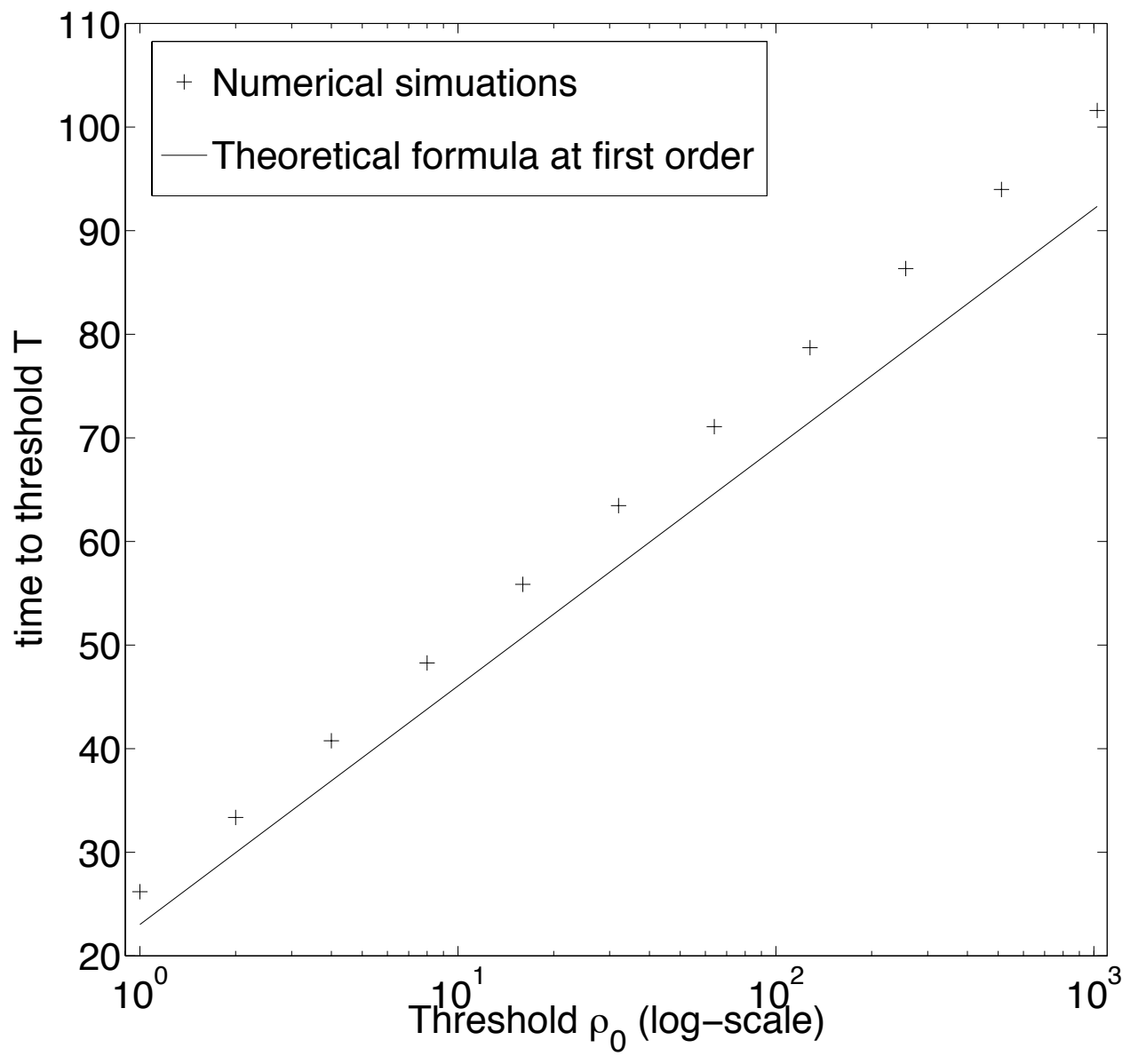

**Fig. 2** Time to threshold $t_{\varrho_0}$ as a function of the threshold parameter $\rho_0$ in log-scale. The crosses are obtained by a direct numerical simulations of the PDE model with constant birth-rate Q. The line is obtained by evaluating formula (27) at first order, i.e. discarding the terms of $o(1/\epsilon)$. Parameters : $b = Q - d = 0.1$, $\epsilon = 0.01$.

*Proof* We introduce $x = e^{t_{\varrho_0}}$ and look for an expansion of the form:

$$x = \frac{x_0}{\epsilon^\beta} + o\left(\frac{1}{\epsilon^\beta}\right) \tag{23}$$

First, we write $x^{b-\Lambda_{\varepsilon,k}} = C_{k,\epsilon}\left(\frac{x_0}{\epsilon^\beta}\right)^{b-\lambda_{0,k}}$ where $C_{k,\epsilon} \to 1$ when $\epsilon \to 0$. Furthermore, using assumption (iii), one controls the convergence of $x^{b-\Lambda_{\varepsilon,k}}$ to 0 uniformly in $k$. Therefore,

$$\frac{\varrho_0}{\epsilon} = \Phi_{\epsilon,0} C_{0,\epsilon}\left(\frac{x}{\epsilon^\beta}\right)^{b-\lambda_{0,0}} - \sum_{k\geq 0}\Phi_{\epsilon,k} + \sum_{k\geq 1}\Phi_{\epsilon,k} x^{b-\Lambda_{\varepsilon,k}} \tag{24}$$

$$= \phi_{0,0}\left(\frac{x_0}{\epsilon^\beta}\right)^{b-\lambda_{0,0}} + O(1) \tag{25}$$

We deduce that

$$\frac{\varrho_0}{\phi_{0,0} x^{b-\lambda_{0,0}}\epsilon} = \frac{1}{\epsilon^{\beta(b-\lambda_{0,0})}} \tag{26}$$

implying that $\beta = 1/(b - \lambda_{0,0})$ and $x_0 = (\varrho_0/\phi_{0,0})^{1/(b-\lambda_{0,0})}$.

In terms of the original parameters of the model, we conclude that the time to produce an output quantity $\varrho_0$ is asymptotically given by:

$$t_{\varrho_0} = \frac{1}{Q_0 - d}\ln\left(\frac{\varrho_0(Q_0 - d)}{< n_I, v_{0,0} > \underline{v}_{0,0}}\frac{1}{\epsilon} + o\left(\frac{1}{\epsilon}\right)\right) \tag{27}$$

This formula relates in a compact form the birth rate $Q_0$, the death rate $d$, the width of the selection function $\epsilon$, the mutation rate $\mu$ and the initial condition $n_I$ to the characteristic time-scale of the B-cell production process. Notice that this formula relies on the assumption that only the first mode grows, meaning that $\mu$ must be larger than $(Q_0 - d)/\pi^2$. Numerical comparison between the time $t_{\varrho_0}$ computed from the numerical solution of the PDE and this formula is displayed in Figure 4.1.4. From the spectral decomposition, we learn that decreasing the mutation rate $\mu$ has the effect of recruiting further modes, whereas a large $\mu$ implies that the evolution of the solution forgets the other modes which were present in the initial condition.

## 4.2 Asymptotic expansion

### *4.2.1 Asymptotic expansion of the solution*

When $Q$ is piecewise constant, using the notation $b = Q_0 - d$, $n(t,x)$, the solution of (1), solves as well the following linear problem, until $\rho(t)$ reaches $\rho_0$:

$$\begin{cases} \partial_t n - \mu \Delta n = (b-s)n & (t,x) \in Q_T, \\ \partial_x n(t,x) = 0 & (t,x) \in (0,T) \times \partial\Omega, \\ n(0,x) = n_I(x) & (t,x) \in \{0\} \times \Omega \end{cases} \tag{28}$$

We look for an approximation of the solution of the system :

$$\begin{cases} \partial_t \mathcal{N} - \mu \Delta \mathcal{N} + s\mathcal{N} = 0, & (t,x) \in Q_T, \\ \partial_x \mathcal{N} = 0, & (t,x) \in (0,T) \times \partial\Omega, \\ \mathcal{N}(0,x) = n_I(x), & (t,x) \in \{0\} \times \Omega. \end{cases} \tag{29}$$

We make the Ansatz :

$$\mathcal{N}_\mu = \mathcal{N}_0\left(t,x,\frac{x}{\varepsilon}\right) + \varepsilon \mathcal{N}_1\left(t,x,\frac{x}{\varepsilon}\right) + \cdots$$

When we plug it in (29), it gives after separating orders of $\varepsilon$ that $\mathcal{N}_0(t,x,y)$ is independent on the fast variable $y$ it solves the homogeneous equations:

$$\begin{cases} \partial_t \mathcal{N}_0 - \text{ă}\mu \Delta \mathcal{N}_0 = 0, & (t,x) \in (0,T) \times \Omega \\ \partial_x \mathcal{N}_0 = 0, & (t,x) \in (0,T) \times \partial\Omega \\ \mathcal{N}_0(0,x) = n_I(x) & (t,x) \in \{0\} \times \Omega, \end{cases}$$

and thus $\mathcal{N}_0(t,x) = \sum_{k\in\mathbb{N}} \overline{n}_I^k \exp(-\lambda_{0,k} t) v_k(x)$.

**Proposition 2** *If $n_I \in L^2(\Omega) \cap L^\infty(\Omega)$ then at zeroth order one can approach $n(t,x)$ solving (28) by $\exp(bt)\mathcal{N}_0(t,x)$ and the error in the $L^\infty((0,T);L^2(\Omega))$ norm is estimated as :*

$$\sup_{t\in(0,T)} \|n(t,\cdot) - \exp(bt)\mathcal{N}_0(t,\cdot)\|_{L^2(\Omega)} \le \varepsilon \exp(bT) T \|n_I\|_{L^\infty(\Omega)}.$$

*Proof* One defines the zero order error :

$$\mathcal{E}_0(t,x) := \mathcal{N}(t,x) - \mathcal{N}_0(t,x)$$

which solves in the strong sense

$$\begin{cases} \partial_t \mathcal{E}_0 - \mu \Delta \mathcal{E}_0 + s\mathcal{E}_0 = -s\mathcal{N}_0, & (t,x) \in (0,T) \times \Omega \\ \partial_x \mathcal{E}_0 = 0, & (t,x) \in (0,T) \times \partial\Omega \\ \mathcal{E}_0(0,x) = 0, & (t,x) \in \{0\} \times \Omega, \end{cases}$$

using standard *a priori* estimates, one obtains that

$$\frac{1}{2}\partial_t \|\mathcal{E}_0(t,\cdot)\|^2_{L^2(\Omega)} \le \|s\mathcal{N}_0\|_{L^2(\Omega)} \|\mathcal{E}_0\|_{L^2(\Omega)}$$

which then by dividing both sides by $\sqrt{\|\mathcal{E}_0\|^2_{L^2(\Omega)} + \delta}$ one gets :

$$\frac{1}{2}\partial_t \sqrt{\|\mathcal{E}_0(t,\cdot)\|^2_{L^2(\Omega)} + \delta} \le \|s\mathcal{N}_0\|_{L^2(\Omega)}$$

which integrated provides after passing to the limit $\delta \to 0$,

$$\|\mathcal{E}_0(t,\cdot)\|_{L^2(\Omega)} \le \int_0^t \left\|s(\cdot)\mathcal{N}_0(\tilde{t},\cdot)\right\|_{L^2(\Omega)} d\tilde{t} \le \varepsilon t \|\mathcal{N}_0\|_{L^\infty((0,t)\times\Omega)} \le \|n_I\|_{L^\infty(\Omega)} \varepsilon t$$

and thus the result follows.

*4.2.2 Estimation of the time to threshold*

Using the above asymptotic expansion, we are now able to estimate the time $t_{\varrho_0}$:

**Theorem 4** *If $s$ satisfies hypotheses 32, and $\mu > b/\pi^2$, then for every given $\varrho_0$ , there exists a time*

$$t_{\varrho_0} := \frac{1}{b} \ln\left(1 + \frac{\varrho_0 b}{\varepsilon \overline{n}_0}\right), \tag{30}$$

*s.t. there exists a constant $C$ independent on $t_{\varrho_0}$ s.t.*

$$|\varrho_0 - \varrho_{\text{out}}(t_{\varrho_0})| \leq |\ln(\varepsilon)|\sqrt{\varepsilon} C(\mu, n_I)$$

*Proof* The previous proposition allows to compute $\varrho_{\text{app}}$, an approximation of $\varrho_{\text{out}}$ which reads :

$$\begin{aligned}\varrho_{\text{app}}(t) &:= \int_0^t \int_\Omega s(x)\mathcal{N}_0(t,x)\exp(bt)dxd\tilde{t} \\ &= \frac{(\exp(bt)-1)}{b}\int_\Omega n_I(x)dx\int_\Omega s(x)dx + \sum_{k\in\mathbb{N}^*}\frac{(\exp((b-\lambda_k)t)-1)}{(b-\lambda_k)} < n_I, v_k >< s, v_k >\end{aligned}$$

which under hypotheses 32 gives

$$\varrho_{\text{app}}(t) := \varepsilon\frac{(\exp(bt)-1)}{b}\overline{n}_I^0 + \sum_{k\in\mathbb{N}^*}\frac{(\exp((b-\lambda_k)t)-1)}{(b-\lambda_k)} < n_I, v_k >< s, v_k >$$

Using Cauchy-Schwartz, one has that

$$\begin{aligned}|\varrho_{\text{out}}(t) - \varrho_{\text{app}}(t)| &\leq \int_0^t\int_\Omega s|n(t,x) - \exp(b\,t)\mathcal{N}_0(\tilde{t},x)|dxd\tilde{t} \leq \int_0^t \|s\|_{L^2(\omega)}\left\|n(t,\cdot) - \exp(b\tilde{t})\mathcal{N}_0(\tilde{t},\cdot)\right\|_{L^2(\Omega)}d\tilde{t} \\ &\leq \varepsilon^{\frac{3}{2}}\|n_I\|_{L^\infty(\Omega)}\int_0^t \exp(b\tilde{t})\tilde{t}d\tilde{t} \leq \varepsilon^{\frac{3}{2}}\|n_I\|_{L^\infty(\Omega)}\frac{t\exp(bt)}{b}\end{aligned}$$

On the other hand using the explicit expression of $\varrho_{\text{app}}$ one writes:

$$\begin{aligned}\left|\varrho_{\text{app}} - \left(\frac{\exp(bt)-1}{b}\right)\overline{n}_I^0\overline{s}\right| &\leq S(t)\sum_{k\in\mathbb{N}^*}|< n, v_k >< s, v_k >| \\ &\leq S(t)\|s\|_{L^2(\Omega)}\|n_I\|_{L^2(\Omega)} \leq S(t)\sqrt{\varepsilon}\|n_I\|_{L^2(\Omega)}\end{aligned} \tag{31}$$

where we denoted $S(t) := \sup_{k\in\mathbb{N}^*}\left|\frac{\exp((b-\lambda_k)t)-1}{(b-\lambda_k)}\right|$. Then two cases occur:

i) either there exists $k_0$ s.t. $\lambda_{k_0} \leq b \leq \lambda_{k_0+1}$ and then

$$S(t) \leq \frac{\exp((b-\mu\pi^2)t)-1}{b-\lambda_{k_0}},$$

ii) or $b < \lambda_1 = \mu\pi^2$ and thus

$$S(t) \leq \frac{2}{(\lambda_1 - b)} = \frac{2}{(\mu\pi^2 - b)}.$$

Using a triangular inequality gives :

$$|\varrho_{\text{out}}(t_{\varrho_0}) - \varrho_0| \leq |\varrho_{\text{out}}(t_{\varrho_0}) - \varrho_{\text{app}}(t_{\varrho_0})| + |\varrho_{\text{app}}(t_{\varrho_0}) - \varrho_0|$$

which, because of the estimates ii) above, gives

$$|\varrho_{\text{out}}(t_{\varrho_0}) - \varrho_0| \leq C\left(\varepsilon^{\frac{3}{2}}t_{\varrho_0}\exp(bt_{\varrho_0}) + \sqrt{\varepsilon}\right) \leq C\sqrt{\varepsilon}|\ln(\varepsilon)|.$$

*Remark 1* Due to (31), one can not, to our knowledge, improve the accuracy when computing $t_{\varrho_0}$ by increasing the order of the asymptotic expansion since the major source of error comes from this step, when integrating in the selection window the zero order term.

*4.2.3 Numerical simulation*

The parameters we choose are $\varrho_0 = 100, b = 2$ together with a random initial condition $n_I$. We compute the error between $t_{\varrho_0}^{h,k}$, the time to reach $\varrho_0$ and the theoretical formula (30). The direct simulation is made using a P1-Finite Element Method algorithm with first order implicit Euler scheme for the time discretization. The results are displayed for various (small) values of $\varepsilon$. We plot, in fig. 3, the error estimates for two values of $\mu$ corresponding to $\mu = b/((k+1/2)^2\pi^2$ for $k \in \{0,1\}$. When $k = 0$ we are in the hypotheses of the latter theorem, whereas for $k = 1$ the theoretical error is not uniform wrt $\varepsilon$. For this specific test-case, the error is comparable in both cases, although greater as $\mu$ becomes smaller as expected. The numerical order of convergence is greater than what is predicted theoretically.

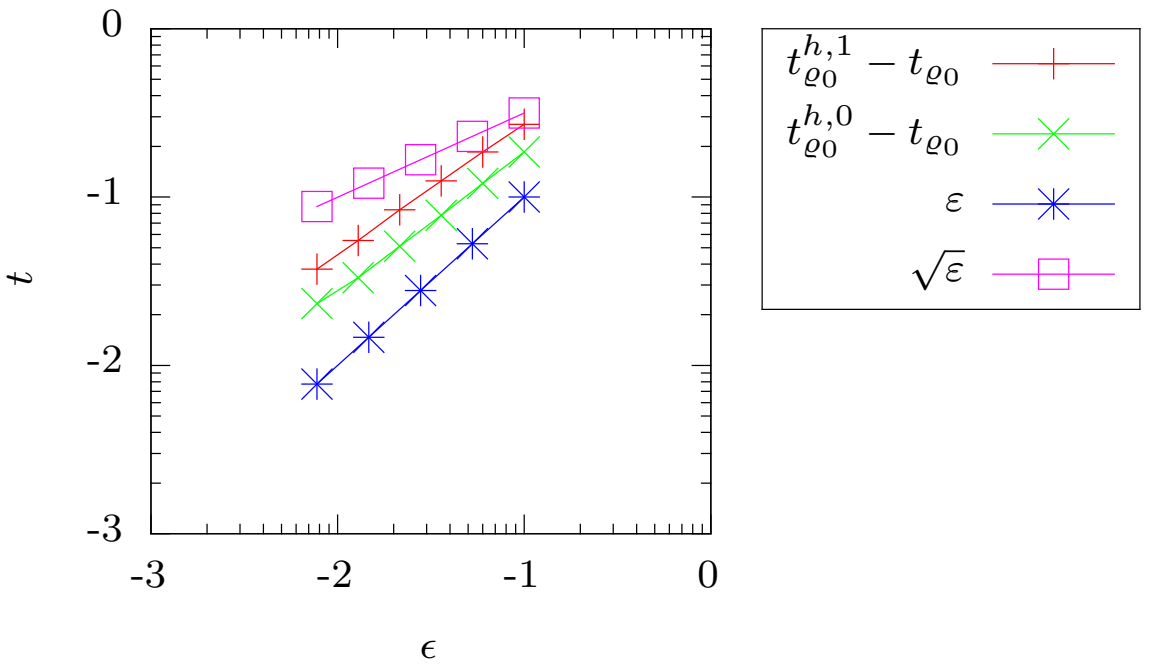

**Fig. 3** Starting from a random initial condition, the difference between the numerical $t_{\varrho_0}^{h,k}$ and the limit time to reach $t_{\varrho_0}$ for various values of $\varepsilon$, the simulations are performed using a P1-FEM from the FEniCS Library [27] with first order implicit Euler scheme for the time discretization.

## 5 Asymptotic analysis for small and large mutation rates

In this section, we focus our attention on the regimes of small and large mutation rates. We first consider the case of initial conditions restricted to a single trait (Dirac initial data), deriving explicit solutions in this case, and providing asymptotics for the time to threshold. Then, we establish asymptotic expansions of the solution in the regimes $\mu \ll 1$ and $\mu \gg 1$.

### 5.1 Explicit solutions and asymptotics for Dirac initial data

*5.1.1 The case of the whole space*

When the size of the domain goes to infinity we face the problem

$$\begin{cases} \partial_t \mathcal{G}_\varepsilon - \mu \Delta \mathcal{G}_\varepsilon + s_\varepsilon \mathcal{G}_\varepsilon = b\mathcal{G}_\varepsilon, & (t,x) \in \mathbb{R}_+ \times \mathbb{R}, \\ \mathcal{G}_\varepsilon(0,x) = \delta_z(x), & x \in (0,1), \end{cases} \tag{32}$$

**Definition 51** *We define a very weak solution of (32) the function $\mathcal{G}_\varepsilon \in L^2(0,T \times \mathbb{R})$ that solves :*

$$\int_{(0,T)\times\mathbb{R}} \mathcal{G}_\varepsilon \left(-\partial_t - \Delta + (s_\varepsilon(x) - b)\right) \varphi dx dt - \varphi(0,z) = 0$$

*for every $\varphi \in C(0,T;H^1(\mathbb{R})) \cap L^2(0,T;H^2(\mathbb{R}))$.*

**Theorem 5** *Provided that $s_\varepsilon \in L^\infty(\mathbb{R})$ there exists a unique very weak function $\mathcal{G}_\varepsilon \in L^2((0,T) \times \mathbb{R})$ solving (32). Moreover one has the comparison principle giving :*

$$\tilde{\mathcal{G}}_{b-s_\infty}(t,x,z) \le \mathcal{G}_\varepsilon(t,x,z) \le \tilde{\mathcal{G}}_b(t,x,z), \;\; a.e\ (t,x) \in (0,T) \times \mathbb{R}$$

*where $\tilde{\mathcal{G}}_b$ is the fundamental solution of the heat equation and reads :*

$$\tilde{\mathcal{G}}_b(t,x,z) = \frac{1}{\sqrt{4\mu\pi t}} \exp\left(bt - \frac{(x-z)^2}{4\mu t}\right),$$

*and $s_\infty := \|s\|_{L^\infty(\mathbb{R})}$.*

The proof relies on duality arguments and the Riesz Theorem and is left to the reader for sake of conciseness.

**Lemma 4** *Setting $J_a(t) := \int_0^t \exp\left(at - \frac{1}{t}\right)\frac{d\tilde{t}}{\sqrt{\tilde{t}}}$ one has*

$$J_a(t) \geq \frac{1}{2e}\exp\left(-\frac{2}{t}\right), \quad \forall t \geq 0$$

*Proof* Using Jensen's inequality one writes

$$\begin{aligned} J_a(t) &\geq \int_0^t \exp\left(-\frac{1}{\tilde{t}}\right)\frac{d\tilde{t}}{\sqrt{\tilde{t}}} = \int_{1/t}^\infty \exp(-z)\frac{dz}{z^{\frac{3}{2}}} \geq \left(\int_{1/t}^\infty \exp(-z)dz\right)\left(\frac{\int_{1/t}^\infty \exp(-z)zdz}{\int_{1/t}^\infty \exp(-z)dz}\right)^{-\frac{3}{2}} \\ &= \exp\left(-\frac{1}{t}\right)\left(\frac{t}{t+1}\right)^{\frac{3}{2}} \geq \exp\left(-\frac{1}{t}\right)\frac{1}{2}\exp\left(-\frac{t+1}{t}\right) = \frac{1}{2e}\exp\left(-\frac{2}{t}\right) \end{aligned}$$

The last estimate comes from estimating $t^{\frac{3}{2}}$ by an exponential, *i.e.* $\forall t \geq 0, \quad \exp\left(-\frac{1}{t}\right) < 2t^{\frac{3}{2}}$.

**Theorem 6** *We suppose that $s_\varepsilon$ is the indicator function of the set $(-\varepsilon,\varepsilon)$ and that $z > \varepsilon$. If we denote the time $t_{\varrho_0}$ s.t. $\rho_\varepsilon(t) := \int_0^t \int_{\mathbb{R}} s(x)\mathcal{G}_\varepsilon(t,x,z)dxdt$ reaches $\varrho_0$, one has that $t_{\varrho_0} \to \infty$ when either $\mu \to 0$ or $\mu \to \infty$. Moreover $t_{\varrho_0}$ does not grow faster than any polynomial wrt $\mu$.*

*Proof* We shall provide a lower bound for $I_l(t) := 2\varepsilon\int_0^t \tilde{\mathcal{G}}_{b-s_\infty}(t,-\varepsilon,z)dt$ and un upper bound for $I_u(t) := 2\varepsilon\int_0^t \tilde{\mathcal{G}}_b(t,\varepsilon,z)dt$ . So that if there exists a time $t_u$ (resp. $t_l$) s.t. $I_l(t_u) \geq \varrho_0$ (resp. $I_u(t_l) \leq \varrho_0$) then $t_{\varrho_0} \leq t_u$ (resp. $t_{\varrho_0} \geq t_l$). We distinguish two cases :

- either $t < t_0 := (z+\varepsilon)/(2\sqrt{(b+s_\infty)\mu})$ and then one can estimate thanks to Lemma 4:

$$I_l(t) \geq \frac{\varepsilon(z+\varepsilon)}{4\mu e\sqrt{\pi}}\exp\left(-\frac{(z+\varepsilon)^2}{2\mu t}\right) =: j_1(t).$$

- or $t > t_0$ and one has

$$\begin{aligned} I_l(t) &\geq j_1(t_0) + \frac{\exp(-(b-s_\infty)t_0)}{\sqrt{\mu\pi}}\int_{t_0}^t \exp((b-s_\infty)\tilde{t})\frac{d\tilde{t}}{\sqrt{\tilde{t}}} \\ &\geq j_1(t_0) + \frac{\exp(-(b-s_\infty)t_0)}{\sqrt{\mu\pi}k!}\int_{t_0}^t (b-s_\infty)^k\tilde{t}^{k-\frac{1}{2}}d\tilde{t} \\ &= j_1(t_0) + \frac{\exp(-(b-s_\infty)t_0)}{\sqrt{\mu\pi}k!(k+\frac{1}{2})}(b-s_\infty)^k\left[\tilde{t}^{k+\frac{1}{2}}\right]_{\tilde{t}=t_0}^{\tilde{t}=t} =: j_2(t). \end{aligned}$$

Inverting this formula one obtains

$$t_u := \begin{cases} j_1^{-1}(\varrho_0) & \text{if } \varrho_0 < j(t_0), \\ j_2^{-1}(\varrho_0) & \text{otherwise }, \end{cases}$$

where

$$\begin{aligned} j_1^{-1}(\varrho_0) &:= \frac{(z+\varepsilon)^2}{-2\ln\left(\frac{\varrho_0 4e\mu}{\varepsilon(z+\varepsilon)}\right)}, \\ j_2^{-1}(\varrho_0) &:= \left(\frac{(k+\frac{1}{2})k!\sqrt{\pi\mu}}{(b-s_\infty)^k}\exp((b-s_\infty)t_0)(\varrho_0 - j_1(t_0)) + t_0^{k+\frac{1}{2}}\right)^{\frac{1}{k+\frac{1}{2}}}, \end{aligned}$$

this proves the second part of the claim. For the lower bound $t_l$, we estimate the time to reach $\varrho_0$ by the upper bound on $I_u$. Indeed we write :

$$I_u(t) \leq \frac{2\varepsilon}{\sqrt{\mu\pi}}\exp\left(bt - \frac{(z-\varepsilon)^2}{4\mu t}\right)\sqrt{t} \leq \frac{2\varepsilon}{\sqrt{\mu\pi}}\exp\left((b+1)t - \frac{(z-\varepsilon)^2}{4\mu t}\right)$$

this gives that

$$t_l = \frac{\ln\omega + \sqrt{\ln^2\omega + (z-\varepsilon)^2(b+1)/\mu}}{2(b+1)}, \quad \omega := \sqrt{\pi\mu}\varrho_0/2\varepsilon$$

which proves the first part of the claim.

In the previous proof the results remain identical if $z < -\varepsilon$ by symmetry. We plot in fig. 4 left, the comparison between our bounds and numerical computations of the integrals $I_l$ and $I_u$ for a given set of data $(b, \varrho_0, \varepsilon, z)$ when $\mu$ varies : $t_{\varrho_0}$ should be comprised between the $t_{I_u}$ and $t_{I_l}$ curves.

#### *5.1.2 The case of a bounded domain*

We look for a solution of the problem : find $\mathcal{G}_\varepsilon$ solution of

$$\begin{cases} \partial_t \mathcal{G}_\varepsilon - \mu\Delta\mathcal{G}_\varepsilon + s_\varepsilon\mathcal{G}_\varepsilon = b\mathcal{G}_\varepsilon, & (t,x) \in \mathbb{R}_+ \times (0,1), \\ \mathcal{G}_\varepsilon'(t,0) = \mathcal{G}_\varepsilon'(t,1) = 0, & t \in \mathbb{R}_+, \\ \mathcal{G}_\varepsilon(0,x) = \delta_z(x), & x \in (0,1), \end{cases} \tag{33}$$

where the support of the dirac mass is located in $z \in (0,1)$.

**Theorem 7** *If $s \in L^\infty(0,1)$ and $b \in \mathbb{R}$, there exists a unique very weak solution $\mathcal{G}_\varepsilon \in L^2((0,T) \times (0,1))$ for every given positive $T$* i.e.

$$\int_{(0,T)\times\mathbb{R}} \mathcal{G}_\varepsilon \left(-\partial_t - \Delta + (s(x) - b)\right)\varphi dx dt - \varphi(0,z) = 0$$

*for every test function $\varphi \in C([0,T]; H^1((0,1))) \cap L^2((0,T); H^2((0,1)))$. Moreover, one has for almost every $(t,x) \in \mathbb{R}_+ \times (0,1)$ that*

$$\tilde{\mathcal{G}}_{b-s_\infty}(t,x,z) \le \mathcal{G}_\varepsilon(t,x,z) \le \tilde{\mathcal{G}}_b(t,x,z)$$

*where upper and the lower bounding functions are fundamental solutions of the heat equation with Neumann boundary conditions :*

$$\begin{cases} \partial_t \tilde{\mathcal{G}}_b - \mu\Delta\tilde{\mathcal{G}}_b = b\tilde{\mathcal{G}}_b, & (t,x) \in \mathbb{R}_+ \times (0,1), \\ \tilde{\mathcal{G}}_b'(t,0) = \tilde{\mathcal{G}}_b'(t,1) = 0, & t \in \mathbb{R}_+, \\ \tilde{\mathcal{G}}_b(0,x) = \delta_z(x), & x \in (0,1), \end{cases} \tag{34}$$

*can be expressed as*

$$\tilde{\mathcal{G}}_b(t,x,z) := \frac{\exp(bt)}{\sqrt{4\pi\mu t}} \sum_{n\in\mathbb{Z}} \left\{ \exp\left(-\frac{(x-2n+z)^2}{4\mu t}\right) + \exp\left(-\frac{(x-2n-z)^2}{4\mu t}\right) \right\}.$$

We need here two technical lemmas :

**Lemma 5** *If $x \in (0,\varepsilon)$ and $z \in (\varepsilon,1)$ and $t \ge 0$ then*

$$\tilde{\mathcal{G}}_b(t,x,z) \le \exp\left(bt - \frac{(z-\varepsilon)^2}{4\mu t}\right)\left(\frac{4}{\sqrt{4\pi\mu t}} + \frac{2}{\sqrt{2(1-\varepsilon)}}\right)$$

*Proof* Thanks to the previous results there exists a $t_1$ s.t. $\tilde{\rho}_b(t_1) > \varrho_0$. This implies that $t_b$, the time that $\tilde{\rho}_b(t) > \varrho_0$ for all $t > t_b$ is greater than $t_1$. As in turn $t_0 \ge t_b$ this ends the proof.

**Lemma 6** *If $(x,z) \in (0,1)^2$, the fundamental solution can be estimated from below as :*

$$\tilde{\mathcal{G}}_b(t,x,z) \ge \exp\left(bt - \frac{(x+z)^2}{4\mu t}\right)\left\{\frac{1}{\sqrt{4\pi\mu t}} + \frac{1}{4\sqrt{3}}\text{erfc}\left(\sqrt{\frac{3}{\mu t}}\right)\right\}.$$

*Proof* We write simply that:

$$\begin{aligned}\sum_{|n|\geq 1} \exp\left(bt - \frac{(x+z-2n)^2}{4\mu t}\right) &\geq \exp\left(bt - \frac{(x+z)^2}{4\mu t}\right) \sum_{|n|\geq 1} \exp\left(-\frac{(x+z)|n|}{\mu t} - \frac{n^2}{\mu t}\right) \\ &\geq 2\exp\left(bt - \frac{(x+z)^2}{4\mu t}\right) \sum_{n\geq 1} \exp\left(-\frac{3n^2}{\mu t}\right) \\ &\geq 2\exp\left(bt - \frac{(x+z)^2}{4\mu t}\right) \int_1^\infty \exp\left(-\frac{3s^2}{\mu t}\right) ds\end{aligned}$$

which ends the proof.

Again we can estimate the time $t_{\varrho_0}$ for which the selected population $\varrho_{\text{out}}$ reaches the threshold value $\varrho_0$.

**Theorem 8** *We suppose that $s_\varepsilon$ is the indicator function of the set $(-\varepsilon, \varepsilon)$ and that $z > \varepsilon$, $b \in \mathbb{R}_+$ s.t. $b > s_\infty := \|s\|_{L^\infty(0,1)}$ and $z > \varepsilon$, then if we denote $t_{\varrho_0}$ the time s.t. $\varrho_{\text{out}}$ reaches the threshold value $\varrho_0$, one has the asymptotics wrt $\mu$ :*

$$\lim_{\mu\to 0} t_{\varrho_0} = +\infty, \quad \lim_{\mu\to\infty} t_{\varrho_0} \in (\underline{t}, \bar{t}),$$

*where the interval $(\underline{t}, \bar{t})$ depends only on the data set $(\varrho_0, \varepsilon, z, b, s_\infty)$, and $0 < \underline{t} < \bar{t} < \infty$.*

*Proof* When $x < \varepsilon < z$ the heat kernel $\tilde{\mathcal{G}}_b$ is monotone wrt $x$. Firstly we compute the lower bound $t_l$ which is provided estimating $I_u(t) := \varepsilon \int_0^t \tilde{\mathcal{G}}_b(t, \varepsilon, z) dt$ from above. By Lemma 5, one has that

$$I_u(t) \leq \varepsilon \frac{\exp\left((b+1)t - \frac{(z-\varepsilon)^2}{4\mu t}\right)}{\sqrt{4\mu\pi}} \left(1 + \sqrt{\frac{2\mu}{(1-\varepsilon)}}\right)$$

and so the time that $I_u$ reaches $\varrho_0$ is then greater than

$$t_l := \frac{\ln\omega + \sqrt{\ln^2\omega + \frac{(z-\varepsilon)^2}{\mu}}}{2(b+1)}, \text{ and } \omega := \frac{2\varrho_0\sqrt{(1-\varepsilon)\mu\pi}}{\varepsilon(\sqrt{1-\varepsilon} + \sqrt{2\mu})}.$$

If $\mu$ is small then

$$t_l \sim \frac{|z-\varepsilon|}{2(b+1)} \frac{1}{\sqrt{\mu}},$$

Whereas if $\mu \to \infty$ ,

$$t_l \to \frac{1}{(b+1)} \ln\left(\frac{\varrho_0\sqrt{2\pi(1-\varepsilon)}}{\varepsilon}\right) =: \underline{t}.$$

For what concerns $t_u$, one has thanks to Lemma 6

$$\tilde{\mathcal{G}}_{b-s_\infty}(t,x,z) \geq \tilde{\mathcal{G}}_{b-s_\infty}(t,0,z) \geq \exp\left((b-s_\infty)t - \frac{z^2}{4\mu t}\right)\left\{\frac{1}{\sqrt{4\pi\mu t}} + \frac{1}{4\sqrt{3}}\operatorname{erfc}\left(\sqrt{\frac{3}{\mu t}}\right)\right\}$$

which allows to write :

$$I_l(t) := \varepsilon \int_0^t \tilde{\mathcal{G}}_{b-s_\infty}(t,0,z) dt \geq I_1(t) + I_2(t),$$

where

$$\begin{aligned} I_1(t) &:= \frac{\varepsilon}{2\sqrt{\pi\mu}} \int_0^{\tilde{t}} \exp\left((b-s_\infty)\tilde{t} - \frac{z^2}{4\mu\tilde{t}}\right) \frac{d\tilde{t}}{\sqrt{\tilde{t}}}, \\ I_2(t) &:= \frac{\varepsilon}{4\sqrt{3}} \int_0^t \exp\left((b-s_\infty)\tilde{t} - \frac{z^2}{4\mu\tilde{t}}\right) \operatorname{erfc}\left(\sqrt{\frac{3}{\mu\tilde{t}}}\right) d\breve{a}\tilde{t}, \end{aligned}$$

and both functions are increasing wrt $t$ which allows us to estimate $I_l$ as

$$I_l(t) \geq I_1(t)\mathbb{1}_{(0,t_0)}(t) + \{I_1(t_0) + I_2(t) - I_2(t_0)\}\,\mathbb{1}_{(t_0,\infty)}(t)$$

where as in the unbounded case $t_0$ represents the time at which the coefficient of the exponential function changes sign *i.e.* $t_0 := z/(2\sqrt{\mu(b-s_\infty)})$. As before when $t < t_0$ the estimate follows the same from Lemma 4. Thus we re-define $j_1(t) := \varepsilon z/(4\mu e)\exp(-z^2/(2\mu t))$. Instead when $t > t_0$, one writes:

$$I_l(t) \geq j_1(t_0) + \frac{\varepsilon}{4\sqrt{3}} \int_{t_0}^{t} \exp((b-s_\infty)(t-t_0))\mathrm{erfc}\left(\sqrt{\frac{3}{\mu \tilde{t}}}\right) d\breve{a}\tilde{t},$$

but because $\mathrm{erfc}(\sqrt{3/(\mu t)})$ is a monotone increasing function wrt $t$ one can again estimate the latter term as :

$$\begin{aligned} I_l(t) &\geq j_1(t_0) + \mathrm{erfc}\left(\sqrt{\frac{3}{\mu t_0}}\right) \frac{\varepsilon}{4\sqrt{3}} \int_{t_0}^{t} \exp((b-s_\infty)(t-t_0)) d\breve{a}\tilde{t} \\ &= j_1(t_0) + \mathrm{erfc}\left(\sqrt{\frac{3}{\mu t_0}}\right) \frac{\varepsilon}{4\sqrt{3}} \frac{(\exp((b-s_\infty)(t-t_0)) - 1)}{(b-s_\infty)} \end{aligned}$$

which leads to the inverse function reading :

$$t_u := \begin{cases} j_1^{-1}(\varrho_0) & \breve{a}\text{if } \varrho_0 < j_1(t_0) \\ t_0 + \frac{1}{(b-s_\infty)} \ln\left(1 + \frac{4\sqrt{3}(b-s_\infty)}{\varepsilon \mathrm{erfc}\left(\sqrt{\frac{3}{\mu t_0}}\right)}(\varrho_0 - j_1(t_0))\right) & \text{otherwise}\breve{a}. \end{cases}$$

when studying the limit of the latter expression when $\mu \to \infty$ one concludes that $j_1(t_0) \to 0$ and $\bar{t} := t_0 + \ln(1 + 4\sqrt{3}(b-s_\infty)\varrho_0/\varepsilon)$.

**Corollary 1** *Under the hypotheses of Theorem 5.1.2, one could give another upper bound for $t_{\rho_0}$ :*

$$t_{\rho_0} \leq t_u^\infty := \begin{cases} j_1^{-1}(\rho_0) & \breve{a}\textit{if } \rho_0 < j_1(t_0) \\ j_2^{-1}(\rho_0) & \textit{otherwise} \end{cases}$$

*where*

$$t_0 = \frac{z}{2\sqrt{b-s_\infty}}, \quad j_1(t_0) := \frac{\varepsilon z}{4\mu e\sqrt{\pi}} \exp\left(-\frac{z\sqrt{b-s_\infty}}{\sqrt{\mu}}\right)$$

*and*

$$\begin{aligned} j_1^{-1}(\varrho_0) &:= \frac{z^2}{-2\ln\left(\frac{\varrho_0 4e\mu}{\varepsilon z}\right)}, \\ j_2^{-1}(\varrho_0) &:= \left(\frac{(k+\frac{1}{2})k!\sqrt{\pi\mu}}{(b-s_\infty)^k} \exp((b-s_\infty)t_0)(\varrho_0 - j_1(t_0)) + t_0^{k+\frac{1}{2}}\right)^{\frac{1}{k+\frac{1}{2}}}, \end{aligned}$$

The proof simply takes into account that the Green function $\mathcal{G}_{b-s_\infty}$ in the bounded case is greater that the Green function in the case of the whole space, applying the same arguments as in section 5.1.1 one concludes. This latter estimates shall improve the upper bound for $\mu$ small.

We plot in fig. 4 right, the comparison between our bounds $t_u$, $t_u^\infty$ and $t_l$ and the direct numerical simulations of $\mathcal{G}_\varepsilon$ solving (33) for a given set of data $(b, \varrho_0, \varepsilon, z)$ when $\mu$ varies.

## 5.2 Asymptotic expansion for large mutation rate

### *5.2.1 The formal result*

When $\mu$, the mutation rate dominates, one sets the decomposition :

$$\mathcal{N}_\mu = \mathcal{N}_0 + \frac{1}{\mu}\mathcal{N}_1 + \cdots$$

We suppose moreover that we focus on solutions for long times so that the time scaling should be

$$\mathcal{N}_\mu(\tilde{t}, x) = n(t, x), \quad (t, x) \in \mathbb{R}_+ \times (0, 1)$$

where $\tilde{t} = \mu t$ and we write the asymptotic expansion of the equations wrt to powers of $1/\mu$ :

– the first order equation :

$$\partial_t \mathcal{N}_0 - \Delta \mathcal{N}_0 = 0$$

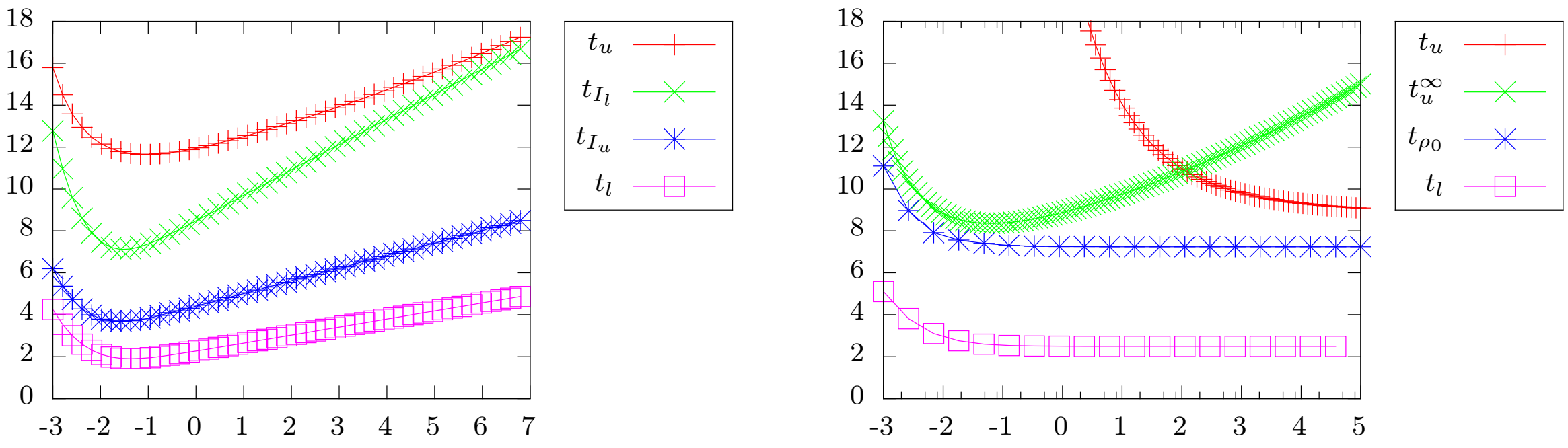

**Fig. 4** The mutation rate $\mu$ is plotted in the log scale on the $x$-axis. The left (resp. right) figure displays the unbounded (resp. bounded) case. On the left side, one shows the upper and lower estimates of $t_{\rho_0}$, the time to reach $\varrho_0$, between these curves we computed numerically $t_{I_l}$ and $t_{I_u}$ in order to validate our estimates. On the right, one plots a numerical computation of $\mathcal{G}_\varepsilon$ solving (33) and the respective theoretical bounds.

– $(j-1)$ $^{\text{th}}$ order :

$$\partial_t \mathcal{N}_j - \Delta \mathcal{N}_j = (b-s)\mathcal{N}_{j-1}, \quad j \in \mathbb{N}^*.$$

Using the spectral decomposition, one gets, for the zeroth order term, that it reads

$$\mathcal{N}_0(t,x) := \sum_{k\in\mathbb{N}} \overline{n}_I^k v_k(x) \exp(-\lambda_{0,k} t), \quad \overline{n}_I^k :=< n_I, v_k >$$

where the brackets denote the scalar product in $L^2((0,1))$, and the eigenvectors are those of the homogeneous problem (12). This expression leads to exponentially decreasing modes and a mean that remains constant.

Denoting $\mathcal{N}_j := \sum_{k\in\mathbb{N}} \gamma_{j,k}(t) v_k(x)$ and $\varGamma_j(t) := (\gamma_{j,k}(t))_{k\in\mathbb{N}}$ one can write the modal equivalent of the equation above :

$$\dot{\varGamma}_j + \mathrm{diag}(\lambda)\varGamma_j = (b - \mathcal{M})\varGamma_{j-1}, \quad (\mathcal{M})_{ik} :=< s v_i, v_k >, \quad \varGamma_j(0) = 0, \ \forall j \geq 1. \tag{35}$$

where $\mathrm{diag}(\lambda) := \mathrm{diag}(\lambda_{0,0}, \dots, \lambda_{0,k}, \dots)$. As the dominant mode is the constant one, we prove below that

$$\gamma_{j,k} \sim \begin{cases} \overline{n}_I^0 \frac{(b-\mathcal{M}_{0,0})^j t^j}{j!} & \text{ăif } k = 0 \\ 0 & \text{ăotherwise} \end{cases} + O(t^{j-1})$$

when $t$ is large. This gives the asymptotic limit :

$$\mathcal{N}_\mu \sim \overline{n}_I^0 \sum_j \frac{(b-\mathcal{M}_{0,0})^j t^j}{j! \mu^j} + R_j t^{j-1} = \overline{n}_I^0 \exp\left( \frac{(b-\mathcal{M}_{0,0})t}{\mu} \right) + \dots .$$

Returning to the original variables and after integration in time this gives that the final formula shall be

$$t_{\varrho_0} := \frac{1}{b - \mathcal{M}_{0,0}} \ln \left( 1 + \frac{\rho_0 \left( b - \mathcal{M}_{0,0} \right)}{\left( \int_0^1 s(x) dx \right) \overline{n}_I^0} \right) + O\left( \frac{1}{\mu} \right).$$

Assuming that $s$ fulfills hypotheses 32 this provides

$$t_{\varrho_0} := \frac{1}{b-\varepsilon} \ln \left( 1 + \frac{\rho_0 \left( b - \varepsilon \right)}{\varepsilon \overline{n}_I^0} \right) + O\left( \frac{1}{\mu} \right).$$

### *5.2.2 Numerical simulations*

We display in fig. 5 for various values of $\mu$, $t_{\varrho_0}$ computed using direct numerical simulations with a random initial data $n_I$, compared with the previous asymptotic values holding for $\mu$ large. In this particular case, the convergence occurs above $\mu \sim O(1)$.

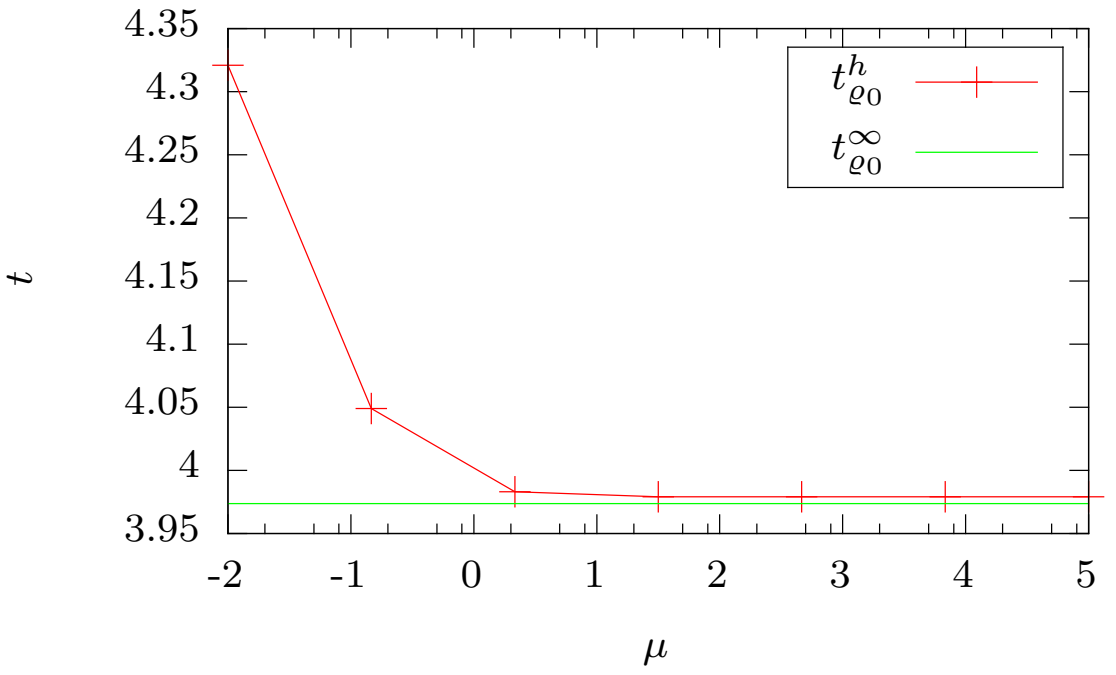

**Fig. 5** Time to reach $\varrho_0$ starting from a random initial condition, for various values of $\mu$, we use a P1-FEM code from the FEniCS Library [27] with first order implicit Euler scheme

*5.2.3 Mathematical proofs*

Setting $\tilde{n}(t,x) := n(t,x)\exp(-(b-\mathcal{M}_{0,0})t)/\overline{n}_I^0$, with $\overline{n}_I^0 := \int_\Omega n_I(x)dx$, it is solution of

$$\begin{cases} \partial_t \tilde{n} - \mu\Delta\tilde{n} = (\mathcal{M}_{0,0} - s)\tilde{n} & (t,x) \in Q_T, \\ \partial_x \tilde{n}(t,x) = 0 & (t,x) \in (0,T)\times\partial\Omega, \\ \tilde{n}(0,x) = n_I(x)/\overline{n}_I^0 & (t,x) \in \{0\}\times\Omega \end{cases} \tag{36}$$

When $\mu$, the mutation rate dominates one sets the decomposition :

$$\mathcal{N}_\mu = \mathcal{N}_0 + \frac{1}{\mu}\mathcal{N}_1 + \cdots$$

We suppose moreover that we focus on solutions for long times so that the time scaling should be

$$\mathcal{N}_\mu(\mu t, x) \sim \tilde{n}(t,x), \quad (t,x) \in \mathbb{R}_+ \times (0,1)$$

and writes the asymptotic expansion of the equations wrt to $1/\mu$ :

- first order equations :

$$\begin{cases} \partial_t \mathcal{N}_0 - \Delta\mathcal{N}_0 = 0, & (t,x) \in (0,T)\times\Omega, \\ \partial_x \mathcal{N}_0(t,x) = 0, & (t,x) \in (0,T)\times\partial\Omega, \\ \mathcal{N}_0(0,x) = n_I(x)/\overline{n}_I^0, & (t,x) \in \{\text{ă}0\}\times\Omega. \end{cases}$$

- zero and higher order equations :

$$\begin{cases} \partial_t \mathcal{N}_j - \Delta\mathcal{N}_j = (\mathcal{M}_{0,0} - s)\mathcal{N}_{j-1}, & (t,x) \in (0,T)\times\Omega, \\ \partial_x \mathcal{N}_j(t,x) = 0, & (t,x) \in (0,T)\times\partial\Omega, \\ \mathcal{N}_j(0,x) = 0, & (t,x) \in \{\text{ă}0\}\times\Omega, \end{cases}$$

where $j \in \mathbb{N}^*$.

**Proposition 3** *For any $n_I \in L^2(\Omega)$, one has*

- *if $j = 0$*

$$\|\mathcal{N}_0(t,\cdot) - 1\|_{L^2(0,1)} \leq c_0 \exp(-\pi^2 t), \quad \forall t \geq 0$$

- *if $j \geq 0$ then*

$$\|\mathcal{N}_j(t,\cdot)\|_{L^2(0,1)} \leq c_j(1+t)^{j-1}, \quad \forall t \geq 0$$

*where the constant $c_j$ does not depend on $t$.*

*Proof* The proof follows in two steps

- If $j = 0$ then the spectral decomposition gives directly the claim.

– For the second part we proceed by induction. If $j=1$, $\mathcal{N}_1$ solves

$$\partial_t \mathcal{N}_1 - \Delta \mathcal{N}_1 = (\mathcal{M}_{0,0} - s)\mathcal{N}_0$$

complemented with homogeneous Neumann boundary and zero initial conditions. Taking the average of this equation provides that $\overline{\mathcal{N}}_1(t) := \int_\Omega \mathcal{N}_1(t,x)dx$ solves

$$\partial_t \overline{\mathcal{N}}_1 = \overline{(\mathcal{M}_{0,0} - s)\mathcal{N}_0} = \overline{(\mathcal{M}_{0,0} - s)\sum_{k\in\mathbb{N}} \gamma_{0,k} v_k} = -\sum_{k\neq 0} \mathcal{M}_{0,k}\gamma_{0,k},$$

but the functions $\gamma_{0,k}$ are explicit and read : $\gamma_{0,k} = \overline{n}_I^k \exp(-\lambda_k t)$. Integrating in time and using Cauchy-Schwartz gives :

$$|\overline{\mathcal{N}}_1|^2 \leq \left( \sup_{k\in\mathbb{N}} |\mathcal{M}_{0,k}\right)^2 | \left| \sum_{k\neq 0ă} \overline{n}_I^k \left( \frac{(1-\exp(-\lambda_k t))}{\lambda_k} \right) \right|^2 \leq c\|s\|^2_{L^\infty(0,1)} \|n_I\|^2_{L^2(0,1)}.$$

For the rest we set $\underline{\mathcal{N}}_1 := \mathcal{N}_1 - \overline{\mathcal{N}}_1$ it solves

$$\partial_t \underline{\mathcal{N}}_1 - \Delta \underline{\mathcal{N}}_1 = (\mathcal{M}_{0,0} - s)\mathcal{N}_0 - \partial_t \overline{\mathcal{N}}_1,$$

which multiplied by $\underline{\mathcal{N}}_1$ and integrated wrt $x$ reads :

$$\frac{1}{2}\partial_t \|\underline{\mathcal{N}}_1\|^2_{L^2(0,1)} + \lambda_1 \|\underline{\mathcal{N}}_1\|^2_{L^2(0,1)} \leq c\|\mathcal{N}_0\|_{L^2(0,1)} \|\underline{\mathcal{N}}_1\|_{L^2(0,1)},$$

by Young's inequality and Gronwall's lemma, the second claim holds for $j=1$.

– we suppose now that for $\ell \leq j-1$ the property is true. Again for the constant mode one has

$$\partial_t \overline{\mathcal{N}}_j = \int_\Omega (\mathcal{M}_{0,0} - s)\mathcal{N}_{j-1} dx \leq \|\mathcal{M}_{0,0} - s\|_{L^2(0,1)} \|\mathcal{N}_{j-1}\|_{L^2(0,1)} \leq ct^{j-2},$$

where the last inequality comes from the induction hypothesis. Integrating the latter inequality in time provides the result for the zero mode. As above, one has

$$\partial_t \underline{\mathcal{N}}_j - \Delta \underline{\mathcal{N}}_j = (\mathcal{M}_{0,0} - s)\mathcal{N}_{j-1} - \partial_t \overline{\mathcal{N}}_j,$$

the latter term being constant in space and $\underline{\mathcal{N}}_j$ being of zero mean value provides after multiplication by $\underline{\mathcal{N}}_j$ and integration wrt $x$ that

$$\frac{1}{2}\partial_t \|\underline{\mathcal{N}}_j\|^2_{L^2(0,1)} + \lambda_1 \|\underline{\mathcal{N}}_j\|^2_{L^2(0,1)} \leq 2\|s\|_{L^\infty(0,1)} \|\mathcal{N}_{j-1}\|_{L^2(0,1)} \|\underline{\mathcal{N}}_j\|_{L^2(0,1)},$$

which again by Young inequality and Gronwall proves the second claim for $\ell = j$. This ends the induction argument.

### *5.2.4 Error estimates*

We define the error as $\mathcal{E}_N := n(t,x)\exp(-bt)/\overline{n}_I^0 - \exp(-\mathcal{M}_{0,0}t)\mathcal{N}_{\mu,N}$ where $\mathcal{N}_{\mu,N} := \sum_{j=0}^{j=N} \frac{1}{\mu^j}\mathcal{N}_j(\mu t, x)$. The error function solves

$$\begin{cases} \partial_t \mathcal{E}_N - \Delta \mathcal{E}_N + s\mathcal{E}_N = \dfrac{(\mathcal{M}_{0,0} - s)}{\mu^N} \mathcal{N}_N(\mu t, x)\exp(-\mathcal{M}_{0,0}t), & (t,x) \in (0,T)\times\Omega \\ \partial_x \mathcal{E}_N = 0, & (t,x) \in (0,T)\times\partial\Omega \\ \mathcal{E}_N(0,x) = 0, & (t,x) \in \{0\}\times\Omega \end{cases}$$

**Theorem 9** *If $n_I \in L^2(\Omega)$ and $s$ satisfies hypotheses 32, one has for any fixed time $T$ that*

$$\|\mathcal{E}_N(T,\cdot)\|_{L^2(\Omega)} \leq \frac{c}{\mu\mathcal{M}_{0,0}^N}.$$

*where the constant $c$ is independent both on $T$ and on $\mu$.*

*Proof* Multiplying the latter equation by $\mathcal{E}_N$ and integrating wrt $x$, one has

$$\partial_t \|\mathcal{E}_N\|^2_{L^2(\Omega)} \leq \frac{2c}{\mu^N} \|\mathcal{N}_N\|_{L^2(\Omega)} \|\mathcal{E}_N\|_{L^2(\Omega)} \exp(-\mathcal{M}_{0,0} t)$$

dividing both sides by $\sqrt{\|\mathcal{E}_N\|^2_{L^2(\Omega)} + \delta}$ one gets :

$$\partial_t \sqrt{\|\mathcal{E}_N\|^2_{L^2(\Omega)} + \delta} \leq \frac{2c}{\mu^N} \frac{\|\mathcal{N}_N\|_{L^2(\Omega)} \|\mathcal{E}_N\|_{L^2(\Omega)}}{\sqrt{\|\mathcal{E}_N\|^2_{L^2(\Omega)} + \delta}} \exp(-\mathcal{M}_{0,0} t) \leq \frac{2c}{\mu^N} \|\mathcal{N}_N(\mu t, \cdot)\|_{L^2(\Omega)} \exp(-\mathcal{M}_{0,0} t)$$

here we use the estimates of Proposition 3 to conclude that

$$\sqrt{\|\mathcal{E}_N(t,\cdot)\|^2_{L^2(\Omega)} + \delta} \leq \sqrt{\delta} + \frac{c}{\mu} \int_0^t s^{N-1} \exp(-\mathcal{M}_{0,0} s) ds \leq \sqrt{\delta} + \frac{c}{\mu \mathcal{M}_{0,0}^N},$$

which finally gives the claim, since the result holds for any arbitrarily small $\delta > 0$.

Turning back to the original function $n(t,x)$ solving (28), the consequence of results above is that

**Corollary 2** *$n_I \in L^2(\Omega)$ and $s$ satisfies hypotheses 32, one has for any fixed time $t$ that*

$$\left\| n(t,\cdot) - \overline{n}_I^0 \exp((b - \mathcal{M}_{0,0})t) \right\|_{L^2(\Omega)} \leq \frac{c |\overline{n}_I^0|}{\mu \mathcal{M}_{0,0}^N} \exp(bt).$$

*where the constant $c$ is independent both on $t$ and on $\mu$.*

*Remark 2* These results show that we understood the limit of $n$ when $\mu$ is large. But the approximation that we built is only first order wrt $\mu$. An interesting and open question is whether one is able to construct a higher order approximation of these asymptotics.

**Theorem 10** *If $n_I \in L^2(\Omega)$ and $s$ satisfies hypotheses 32, for any $\eta > 0$, there exists a $\mu$ big enough, s.t. if*

$$t_{\varrho_0} := \frac{1}{b - \mathcal{M}_{0,0}} \ln \left( 1 + \frac{\rho_0 \left(b - \mathcal{M}_{0,0}\right)}{\left(\int_0^1 s(x) dx\right) \overline{n}_I^0} \right),$$

*then $|\varrho_{\text{out}}(t_{\varrho_0}) - \varrho_0| \leq \eta$.*

## 5.3 Asymptotic expansion for small mutation rates

At the contrary, if we consider $\mu$ small, one can decompose the solution of (28)

$$n(t,x) \sim \sum_{j \in \mathbb{N}} \mu^j \mathcal{N}_j(t,x)$$

where the different terms $\mathcal{N}_j$ solve :

– if $j = 0$

$$\begin{cases} \partial_t \mathcal{N}_0 = (b - s(x)) \mathcal{N}_0, \quad (t,x) \in (0,T) \times \Omega, \\ \mathcal{N}_0(0,x) = n_I(x). \end{cases} \tag{37}$$

– if $j \geq 1$

$$\begin{cases} \partial_t \mathcal{N}_j = (b - s(x)) \mathcal{N}_j + \Delta \mathcal{N}_{j-1}, \quad (t,x) \in (0,T) \times \Omega, \\ \mathcal{N}_j(0,x) = 0. \end{cases} \tag{38}$$

### 5.3.1 Formal computations

One can solve explicitly $\mathcal{N}_0$, which reads :

$$\mathcal{N}_0(t,x) = \exp((b - s(x))t) n_I(x).$$

We compute a zero order approximation of $\varrho_{\text{out}}$ :

$$\varrho_{\text{out}} \sim \int_0^t \int_\Omega s(x) \exp((b - s(x))\tilde{t}) n_I(x) dx d\tilde{t} =: \rho_{\text{app}}(t).$$

For the particular case when $s$ satisfies hypotheses 32 one recovers

$$\rho_{\text{app}}(t) = \frac{\exp((b-1)t) - 1}{(b-1)} \int_\Omega s(x) n_I(x) dx,$$

which gives then an explicit formula for $\mu$ small

$$t_{\varrho_0} := \frac{1}{(b-1)} \log\left(1 + \frac{(b-1)\varrho_0}{\int_\Omega s(x) n_I(x) dx}\right). \tag{39}$$

### 5.3.2 Numerical simulations

We display in fig. 6 $t_{\varrho_0}$, the time to reach $\varrho_0$ using a P1-FEM code with first order implicit Euler scheme for the time discretization, for various values of $\mu$ and for a given random initial data $n_I \in L^2(\Omega)$. We observe that the convergence towards the limit value $t_{\varrho_0}^{\mu=0}$ provided by our latter formula is rather slow.

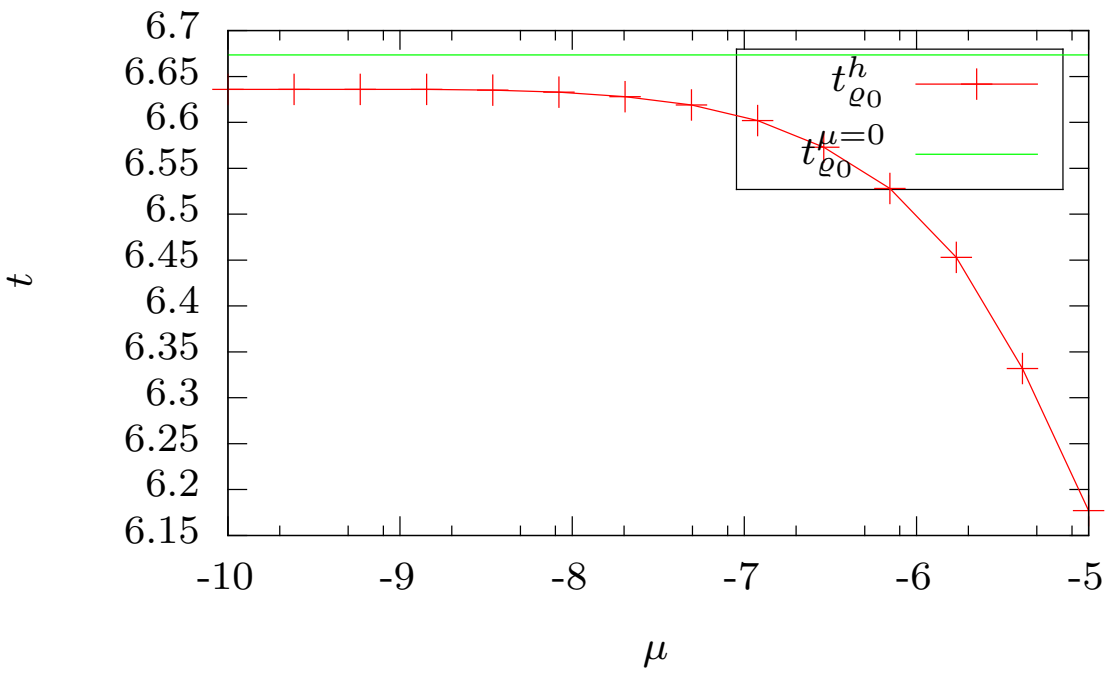

**Fig. 6** Time to reach $\varrho_0$ starting from a random initial condition, for various values of $\mu$, we use again the same code as in fig. 5.

### 5.3.3 Rigorous proofs

**Lemma 7** *We suppose that $s$ satisfies the hypotheses 32. For any given $\Psi \in L^2((0,T) \times \Omega)$ there exists a unique solution $\varphi \in C((0,T); H^1(\Omega)) \cap L^2((0,T); H^2(\Omega))$ solving*

$$\begin{cases} \partial_t \varphi - \mu \Delta \varphi + s\varphi = \Psi(t,x), & (t,x) \in Q_T \\ \partial_x \varphi(t,x) = 0, & (t,x) \in (0,T) \times \partial\Omega \\ \varphi(0,x) = 0, & (t,x) \in \{ǎ0\} \times \Omega, \end{cases} \tag{40}$$

*and one has continuity wrt the $L^2((0,T); H^2(\Omega))$ norm :*

$$\|\varphi\|_{L^2((0,T);H^2(\Omega))} \leq C \|\Psi\|_{L^2(Q_T)}.$$

The proof, based on the Galerkin decomposition, is classical and can be found in [21] Chap III, sect. 6, p.172-178.

**Theorem 11** *If $n_I \in L^2(\Omega)$ and $s$ satisfies hypotheses 32, then one has an $L^2(Q_T)$ error estimates :*

$$\|n - \mathcal{N}_0\|_{L^2(Q_T)} \leq \mu \sqrt{\frac{\exp(2bT) - 1}{2b}} \|n_I\|_{L^2(\Omega)}$$

*Proof* We rescale the problem (28) and the asymptotic expansion so to drop the term containing $b$ in the respective equations. Then one rescales back at the end of the proof. Let define $Y := C^0((0,T); H^1(\Omega)) \cap L^2((0,T); H^2(\Omega))$. Considering (37), the equation is satisfied also in the weak sense namely, for any $\varphi \in Y$ one has

$$\left[ \int_\Omega \mathcal{N}_0(t,x)\varphi(t,x)dx \right]_{t=0}^{t=T} - \int_{Q_T} \mathcal{N}_0(t,x) \left( \partial_t \varphi - s(x)\varphi \right) dxdt = 0.$$

On the other hand, one has also that $n$ solving (28) satisfies as well for any $\varphi \in Y$

$$\left[ \int_\Omega n(t,x)\varphi(t,x)dx \right]_{t=0}^{t=T} - \int_{Q_T} n(t,x) \left( \partial_t \varphi + \mu \Delta \varphi - s(x)\varphi \right) dxdt = 0.$$

This gives when setting $\mathcal{E}(t,x) := n(t,x) - \mathcal{N}_0(t,x)$ that it satisfies for any $\varphi \in Y$

$$\left[ \int_\Omega \mathcal{E}(t,x)\varphi(t,x)dx \right]_{t=0}^{t=T} - \int_{Q_T} \mathcal{E}(t,x) \left( \partial_t \varphi + \mu \Delta \varphi - s(x)\varphi \right) dxdt = \mu \int_{Q_T} \mathcal{N}_0(t,x) \Delta \varphi(t,x) dxdt.$$

For any $\Psi \in L^2((0,T) \times \Omega)$, there exists $\varphi \in Y$, the forward form solving (40) and insert the backward expression of $\varphi$ in the latter weak form which gives :

$$\int_{Q_T} \mathcal{E}(t,x)\Psi(t,x)dxdt \leq \mu \|\mathcal{N}_0\|_{L^2(Q_T)} \|\Delta\varphi\|_{L^2(Q_T)} \leq C\mu \|\mathcal{N}_0\|_{L^2(Q_T)} \|\Psi\|_{L^2(Q_T)},$$

which holds for any $\Psi \in L^2(Q_T)$. Taking the supremum over all functions in this latter space ends the proof.

**Corollary 3** *Under the same hypotheses as in Theorem (11), the time $t_{\varrho_0}$ being defined in (39), one has the error estimates : there exists a constant $C(s, n_I, b) > 0$, independent of $\mu$, s.t.*

$$|\varrho_{\text{out}}(t_{\varrho_0}) - \varrho_0| \leq \mu C(s, n_I, b) \|n_I\|_{L^2(\Omega)} \|s\|_{L^2(\Omega)}.$$

The proof comes easily combining a triangular inequality as in the proof of Theorem 4, results above and the definition of $t_{\varrho_0}$. As the latter time does not depend on $\mu$, the claim follows straightforwardly.